\documentclass[12pt,a4paper,reqno]{article}
\usepackage[utf8]{inputenc}
\usepackage[T1]{fontenc}
\usepackage{amsmath}
\usepackage{graphicx}
\usepackage{amssymb}
\usepackage{mathrsfs}  
\usepackage{amsfonts}
\usepackage{amsthm}
\usepackage[left=2cm, right=2cm, top=3cm, bottom=3.5cm]{geometry}
\usepackage[colorlinks=true,linkcolor=cyan,citecolor=magenta,bookmarks=true]{hyperref}

\DeclareMathAlphabet{\mathcal}{OMS}{cmsy}{m}{n}

\makeatletter

\DeclareFontFamily{OMX}{MnSymbolE}{}
\DeclareSymbolFont{MnLargeSymbols}{OMX}{MnSymbolE}{m}{n}
\SetSymbolFont{MnLargeSymbols}{bold}{OMX}{MnSymbolE}{b}{n}
\DeclareFontShape{OMX}{MnSymbolE}{m}{n}{
    <-6>  MnSymbolE5
   <6-7>  MnSymbolE6
   <7-8>  MnSymbolE7
   <8-9>  MnSymbolE8
   <9-10> MnSymbolE9
  <10-12> MnSymbolE10
  <12->   MnSymbolE12
}{}
\DeclareFontShape{OMX}{MnSymbolE}{b}{n}{
    <-6>  MnSymbolE-Bold5
   <6-7>  MnSymbolE-Bold6
   <7-8>  MnSymbolE-Bold7
   <8-9>  MnSymbolE-Bold8
   <9-10> MnSymbolE-Bold9
  <10-12> MnSymbolE-Bold10
  <12->   MnSymbolE-Bold12
}{}

\let\llangle\@undefined
\let\rrangle\@undefined
\DeclareMathDelimiter{\llangle}{\mathopen}%
                     {MnLargeSymbols}{'164}{MnLargeSymbols}{'164}
\DeclareMathDelimiter{\rrangle}{\mathclose}%
                     {MnLargeSymbols}{'171}{MnLargeSymbols}{'171}

\makeatother

\theoremstyle{thmstyleone}
\newtheorem{theorem}{Theorem}[subsection]
\newtheorem{lemma}[theorem]{Lemma}
\newtheorem{proposition}[theorem]{Proposition}
\newtheorem{corollary}[theorem]{Corollary}

\theoremstyle{thmstylethree}
\newtheorem{definition}[theorem]{Definition}

\theoremstyle{thmstyletwo}
\newtheorem{remark}[theorem]{Remark}

\numberwithin{equation}{subsection}

\newcommand{\build}[3]{\mathrel{\mathop{\kern 0pt#1}\limits_{#2}^{#3}}}

\def\geq{\geqslant}
\def\leq{\leqslant}

\def\epsilon{\varepsilon}

\def\bar{\overline}

\begin{document}

\author{Thibaut Lemoine\thanks{Universit\'e de Strasbourg, CNRS, UMR 7501 -- Institut de Recherche Math\'ematique Avanc\'ee, 7 rue Ren\'e Descartes, 67000 Strasbourg, France. thibaut.lemoine@math.unistra.fr}}

\title{Determinantal point processes on complex manifolds: Construction and limit theorems}

\maketitle

\abstract{We develop a coordinate-free probabilistic framework for determinantal point processes associated with Bergman kernels on compact complex manifolds. The basic issue is that Bergman kernels are naturally line-bundle-valued: $B_k(x,y)\in\operatorname{Hom}(L_y^k,L_x^k)$. Hence the usual determinantal formula for correlation functions is not literally a scalar determinant unless one first gives it an intrinsic meaning. We rigorously define this determinant and prove that every finite-dimensional Hilbert space of sections of a Hermitian line bundle gives rise to a genuine finite-rank projection determinantal point process on the base manifold. We then isolate a collection of finite-dimensional transfer principles showing how diagonal asymptotics, near-diagonal asymptotics, Schur complements, Toeplitz trace expansions and determinant asymptotics are converted into probabilistic statements. Specializing to $H^0(M,L^k)$, this gives the Bergman ensemble as the geometric analogue of an orthogonal polynomial ensemble, and some of the transfer principles allow us to recover previously known results of Berman.}

\tableofcontents

\section{Introduction}

Determinantal point processes (DPPs) are random point configurations whose correlation functions are given by determinants.  They were introduced by Macchi \cite{Mac2} and have become central in probability, random matrix theory and statistical mechanics; see for instance \cite{Sos,Lyo,ShiTak,Joh,HKPV,DPP-Fermion}.  A particularly important and elementary source of DPPs is the following finite-dimensional projection construction.  If $H\subset L^2(E,\mu)$ is an $N$-dimensional Hilbert space of functions on some Polish space with orthonormal basis $\psi_1,\ldots,\psi_N$, then the density
\[
\frac1{N!}
\left|\det\bigl(\psi_i(x_j)\bigr)_{1\leq i,j\leq N}\right|^2
\,d\mu(x_1)\cdots d\mu(x_N)
\]
defines a DPP with projection kernel $K_H(x,y)=\sum_{\ell=1}^N\psi_\ell(x)\overline{\psi_\ell(y)}.$ This simple construction underlies, among many examples, orthogonal polynomial ensembles and Christoffel--Darboux kernels in random matrix theory \cite{Meh,For,Dei,Joh,Bor,BS}.

The purpose of this paper is to formulate this projection-DPP construction intrinsically in the setting of complex geometry.  Let $M$ be a compact complex manifold and let $L\to M$ be a Hermitian holomorphic line bundle. For the finite-dimensional space $H_k=H^0(M,L^k),$ the orthogonal projection onto $H_k$ is represented by the Bergman kernel. However, unlike the scalar kernels of the usual DPP formalism, the Bergman kernel is canonically a line-bundle-valued object, thus the expression $\det\bigl(B_k(x_i,x_j)\bigr)_{1\leq i,j\leq m}$ does not have an immediate scalar meaning.  The first point of the paper is that, because the fibers of $L$ are one-dimensional, this expression has a canonical coordinate-free interpretation: it is the determinant of the endomorphism of $L^k_{x_1}\oplus\cdots\oplus L^k_{x_m}$ whose $(i,j)$-block is $B_k(x_i,x_j)$.  This intrinsic determinant is the scalar object which replaces the usual determinant of a scalar kernel.

The central construction result is Theorem~\ref{thm:intrinsic-determinantal-property}.  It states that if $H$ is any finite-dimensional Hilbert space of measurable sections of a Hermitian line bundle $L\to M$, then the squared norm of the associated Slater section defines a probability measure on $M^N$, and the corresponding point process is a genuine projection determinantal point process on $M$, in the sense that its correlation functions are scalar functions given by the intrinsic determinant of the line-bundle-valued projection kernel
\[
B_H(x,y)\in \operatorname{Hom}(L_y,L_x).
\]
Hence, no choice of local frame, gauge or trivialization enters the definition of the process.  Local scalar kernels may be used for computations, but the correlation functions themselves are the coordinate-free determinants of endomorphisms of
$L_{x_1}\oplus\cdots\oplus L_{x_m}$.

The second part of the paper turns this construction into a probabilistic framework.  The finite-dimensional identities which make projection DPPs useful remain valid in the intrinsic setting: linear statistics are controlled by compressed multiplication operators, multiplicative functionals by finite-dimensional determinants, and reduced Palm measures by Schur complements.  In the Bergman case, the Palm identity has a particularly simple geometric meaning: conditioning on the presence of points $p_1,\ldots,p_\ell$ replaces the one-particle space $H^0(M,L^k)$ by the subspace of sections vanishing at those points, thus Palm conditioning is exactly the operation of imposing zeros on holomorphic sections.

The third part of the framework is the collection of transfer principles proved in Section~\ref{sec:black-box-transfer}.  These results isolate the finite-dimensional probabilistic mechanism which converts asymptotic information on projection kernels into limit theorems for point processes. Theorem~\ref{thm:diagonal-transfer} shows that convergence of the normalized one-point intensity implies convergence in probability of empirical measures.  Theorem~\ref{thm:local-kernel-transfer} transfers local kernel asymptotics to local correlation asymptotics, and Proposition~\ref{prop:local-process-transfer} gives a corresponding criterion for convergence of the rescaled point processes. Proposition~\ref{prop:palm-transfer} gives the Palm analogue by applying Schur complements to the local kernels.  Finally, Proposition~\ref{prop:cumulant-trace-formula} expresses cumulants of linear statistics through traces of compressed multiplication operators, while Theorem~\ref{thm:ldp-transfer} gives an abstract large-deviation transfer theorem from determinant asymptotics. Although each transfer mechanism is already known from the viewpoint of finite-rank DPPs, we reformulate them in a form adapted to intrinsic line-bundle-valued kernels. This makes explicit which analytic input is needed for each probabilistic conclusion.

Applied to Bergman ensembles, these abstract results produce the following dictionary:
\[
\begin{array}{c|c}
\text{analytic input} & \text{probabilistic output} \\
\hline
\text{diagonal Bergman expansion} & \text{law of large numbers for empirical measures} \\
\text{near-diagonal expansion} & \text{local correlation and local fluctuation expansions} \\
\text{off-diagonal localization} & \text{variance bounds for linear statistics} \\
\text{Schur-complement asymptotics} & \text{Palm limits} \\
\text{Toeplitz trace expansions} & \text{cumulant expansions} \\
\text{Gram determinant asymptotics} & \text{large deviations}
\end{array}
\]
The analytic inputs in the left column are results from Bergman kernel theory, Toeplitz quantization and pluripotential theory.  The probabilistic part of the paper is the finite-dimensional mechanism which turns those inputs into the conclusions in the right column.

Some of the Bergman-ensemble consequences obtained in this way are already known, often under weaker analytic assumptions, from Berman's work on determinantal processes and fermions on complex manifolds \cite{Ber6,Ber7}. The aim here is not to prove new Bergman-kernel asymptotics. Instead, the paper makes explicit the intrinsic projection-DPP framework underlying these results and formulates a modular interface between complex-geometric kernel asymptotics and probabilistic limit theorems. In particular, any improvement of the analytic input can be inserted into the same transfer mechanism. For instance, recent asymptotics for partial Bergman kernels and their associated determinantal processes suggest that the same interface can be useful beyond the full Bergman projection; see \cite{Ioo25}.

This intrinsic line-bundle-valued setting should be distinguished from the scalar Bergman-kernel DPPs associated with Hilbert spaces of holomorphic functions on domains, such as the processes studied by Bufetov, Fan and Qiu \cite{BQ,BFQ,Buf23}.  In those works the Bergman kernel is a scalar kernel on the underlying space.  Here the natural object is instead a kernel with values in $\operatorname{Hom}(L_y^k,L_x^k)$, and the scalar correlation functions arise only after applying the intrinsic determinant. The construction is also natural from the point of view of fermionic many-body states and the integer quantum Hall effect \cite{DK,Kle,Kle2}.  In that language, $H^0(M,L^k)$ is the lowest Landau level, the Slater section is the filled fermionic state, and the Bergman ensemble is the position process of the fermions.  The line-bundle-valued nature of the Bergman kernel is then gauge covariance, while the intrinsic determinant gives the gauge-invariant correlation density.

The paper is organized as follows.  Section~\ref{sec:projection-DPPs} recalls the finite-rank projection-DPP identities used throughout the paper and extends their construction to Hilbert spaces of sections of a Hermitian line bundle. Section~\ref{sec:black-box-transfer} collects a number of standard finite-dimensional projection-DPP identities and packages them as transfer principles, so that they can be used as a modular interface between line-bundle-valued Bergman kernel asymptotics and probabilistic limit theorems.  Section~\ref{sec:Bergman-DPP} applies this framework to Bergman ensembles associated with $H^0(M,L^k)$.  The analytic inputs in the final section are quoted from Bergman kernel theory, Toeplitz calculus and pluripotential theory; the point is the way in which they enter the general probabilistic mechanism.

\section{Projection DPPs with line-bundle-valued kernels}
\label{sec:projection-DPPs}

In this section we recall the basic definitions of determinantal point processes, with a special focus on a particular subclass: finite-rank projection determinantal point processes. This is the probabilistic structure behind orthogonal polynomial ensembles, the eigenvalues of many random matrix models, as well as the Bergman ensembles that we will define later.

\subsection{Determinantal point processes and reproducing kernels}

In this paragraph we recall the usual definition and properties of projection DPPs, without proofs. We refer to references such as \cite{Sos,Lyo,ShiTak,Joh,HKPV} for further details.

Let $E$ be a locally compact Polish space and let $\mu$ be a Radon measure on $E$. We denote by $\operatorname{Conf}(E)$ the set of locally finite subsets of $E$. A simple point process on $E$ is a random variable $X : \Omega \longrightarrow \operatorname{Conf}(E)$, defined on a probability space $(\Omega,\mathcal{F},\mathbb P)$. One may also define $X$ as the random counting measure $\mathcal X=\sum_{x\in X}\delta_x.$ In the finite-particle case, which will be the main case in this paper, we write $\mathcal X=\sum_{i=1}^N \delta_{X_i},$ where $X_1,\ldots,X_N$ are random points of $E$. Since the configurations are unordered, the law of $(X_1,\ldots,X_N)$ will always be symmetric.

For $m\geq 1$, the $m$-th factorial moment measure is the measure on $E^m$ given by
\[
    \alpha_m(F)=\mathbb E\left[\sum_{\substack{x_1,\ldots,x_m\in X\\ x_i\neq x_j}}F(x_1,\ldots,x_m)
    \right],
\]
for every non-negative measurable function $F$ on $E^m$. If $\alpha_m$ is absolutely continuous with respect to $\mu^{\otimes m}$, its density $\rho_m$ is called the $m$-point correlation function:
\[
    \mathbb E\left[\sum_{\substack{x_1,\ldots,x_m\in X\\ x_i\neq x_j}}F(x_1,\ldots,x_m)\right]=\int_{E^m}F(x_1,\ldots,x_m)\rho_m(x_1,\ldots,x_m)\,d\mu(x_1)\cdots d\mu(x_m).
\]

If $X$ has exactly $N$ points almost surely and if $(X_1,\ldots,X_N)$ has symmetric density $p_N$ with respect to $\mu^{\otimes N}$, then
\[
    \rho_m(x_1,\ldots,x_m)=\frac{N!}{(N-m)!}\int_{E^{N-m}}p_N(x_1,\ldots,x_N)\,d\mu(x_{m+1})\cdots d\mu(x_N)
\]
for $m\leq N$, and $\rho_m=0$ for $m>N$.

\begin{definition}
Let $K:E\times E\to \mathbb C$ be a measurable kernel. A simple point process $X$ on $E$ is called a \emph{determinantal point process} with kernel $K$, with respect to $\mu$, if its correlation functions satisfy
\begin{equation}
\rho_m(x_1,\ldots,x_m)=\det\big(K(x_i,x_j)\big)_{1\leq i,j\leq m}
\end{equation}
for every $m\geq 1$.
\end{definition}
We shall often abbreviate determinantal point process as DPP. In full generality, the existence of a DPP with a prescribed kernel requires spectral assumptions on the integral operator associated with $K$ \cite{Sos}. In this paper, however, we shall mainly use the following finite-dimensional situation, where the construction is completely explicit. Let $H$ be an $N$-dimensional subspace of $L^2(E,\mu)$, and let $\psi_1,\ldots,\psi_N$ be an orthonormal basis of $H$. The orthogonal projection $\Pi_H:L^2(E,\mu)\longrightarrow H$ has integral kernel
\[
    K_H(x,y)
    =
    \sum_{\ell=1}^N \psi_\ell(x)\overline{\psi_\ell(y)},
\]
which is Hermitian, reproducing on $H$, and does not depend on the choice of the orthonormal basis.

\begin{definition}\label{def:projection-DPP}
Let $H\subset L^2(E,\mu)$ be an $N$-dimensional subspace with projection kernel $K_H$. The \emph{projection determinantal point process} associated with $H$ is the probability measure on $E^N$ with density
\[
    p_H(x_1,\ldots,x_N)=\frac{1}{N!}\left|\det\big(\psi_i(x_j)\big)_{1\leq i,j\leq N}\right|^2
\]
with respect to $\mu^{\otimes N}$, where $\psi_1,\ldots,\psi_N$ is any orthonormal basis of $H$.
\end{definition}

The normalization follows from Andreief's identity, also known as the generalized Cauchy--Binet identity \cite[Prop. 2.10]{Joh}:
\[
    \int_{E^N}\left|\det\big(\psi_i(x_j)\big)_{1\leq i,j\leq N}\right|^2\,d\mu(x_1)\cdots d\mu(x_N)=N!.
\]
Thus $p_H$ is indeed a probability density. If $E=\mathbb R$ or $E=\mathbb C$, and if $H$ is generated by the first $N$ orthonormal polynomials associated with a measure $\mu$, then $K_H$ is the usual Christoffel--Darboux kernel \cite{Dei}. The corresponding projection DPP is the classical orthogonal polynomial ensemble, and the fact that it is indeed a DPP can be found for instance in \cite{Bor,BS}. A large part of the probabilistic information of a projection DPP is already encoded in simple operator identities. Let $f:E\to\mathbb C$ be a bounded measurable function and define the linear statistic $S_f=\sum_{x\in X} f(x).$ Let $M_f$ denote the multiplication operator by $f$ on $L^2(E,\mu)$, and let $T_f=\Pi_H M_f \Pi_H$ be its compression to $H$. Equivalently, $T_f$ is the $N\times N$ matrix with coefficients
\[
    (T_f)_{ij}=\int_E f(x)\psi_i(x)\overline{\psi_j(x)}\,d\mu(x).
\]

\begin{proposition}\label{prop:linear-statistics-projection-DPP}
For every bounded measurable function $f$,
\begin{equation}
\mathbb E[S_f]=\operatorname{Tr}(T_f)=\int_E f(x)K_H(x,x)\,d\mu(x).
\end{equation}
If $f$ is real-valued, then
\begin{equation}\label{eq:variance-dpp}
\operatorname{Var}(S_f)=\frac12\int_{E^2}|f(x)-f(y)|^2|K_H(x,y)|^2\,d\mu(x)d\mu(y).
\end{equation}
\end{proposition}



The same operator $T_f=\Pi_HM_f\Pi_H$ also controls multiplicative functionals: if $g:E\to\mathbb C$ is bounded and $1+g$ is non-negative, then
\begin{equation}\label{eq:mult-functional}
\mathbb E\left[\prod_{x\in X}(1+g(x))\right]=\det_H\big(I_H+T_g\big),
\end{equation}
where $\det_H$ denotes the determinant on the finite-dimensional space $H$. In particular, for a bounded real function $f$,
\begin{equation}\label{eq:exp-lin-stat}
\mathbb E\left[e^{S_f}\right]=\det_H\big(I_H+T_{e^f-1}\big).
\end{equation}
This formula will later be useful when comparing changes of weights in Bergman ensembles with determinants of Gram or Toeplitz matrices.

Projection DPPs are stable under conditioning. This fact will acquire a geometric interpretation later: for Bergman DPPs, conditioning on the presence of points corresponds to restricting holomorphic sections to vanish at those points. Let $p_1,\ldots,p_m\in E$ be distinct points such that $\det\big(K_H(p_a,p_b)\big)_{1\leq a,b\leq m}\neq 0.$ Define the $m\times m$ Gram matrix $G_{\mathbf p}=\big(K_H(p_a,p_b)\big)_{1\leq a,b\leq m}.$ The reduced Palm kernel at $\mathbf p=(p_1,\ldots,p_m)$ is
\[
    K_H^{!\mathbf p}(x,y)=K_H(x,y)-\sum_{a,b=1}^mK_H(x,p_a)(G_{\mathbf p}^{-1})_{ab}K_H(p_b,y).
\]
This is the Schur complement of $G_{\mathbf p}$ in the kernel $K_H$.

\begin{proposition}\label{prop:palm-projection-DPP}
The reduced Palm process of the projection DPP associated with $H$, conditioned to contain $p_1,\ldots,p_m$ and with these points removed, is again a projection DPP, with kernel $K_H^{!\mathbf p}$.
\end{proposition}

\subsection{Line-bundle-valued kernels and intrinsic determinants}
\label{sec:line-bundle-kernels}

In the previous section, projection DPPs were described using scalar kernels $K(x,y)$. However, Bergman kernels on complex manifolds are slightly different: if $L\to M$ is a Hermitian line bundle, the Bergman kernel is not naturally a complex-valued function on $M\times M$. Rather, it is a kernel with values in linear maps between fibers: $B(x,y)\in \operatorname{Hom}(L_y,L_x).$ This creates a small but important point: expressions such as
\[
    \det\big(B(x_i,x_j)\big)_{1\leq i,j\leq m}
\]
are not ordinary scalar determinants and might seem ill-defined at first glance. Nevertheless, because $L$ has one-dimensional fibers, there is a canonical way to interpret this determinant as a scalar. The point of the present section is to extend the existing theory of scalar projection DPPs to DPPs whose kernels are related to line bundles.

Throughout this section, $M$ is a complex manifold endowed with a Radon measure $\mu$, and $L\to M$ is a Hermitian complex line bundle. We write $h_x$ for the Hermitian product on the fiber $L_x$, and we choose the convention that Hermitian products are linear in the first variable. Let $H$ be an $N$-dimensional Hilbert space of measurable sections of $L$, endowed with the inner product inherited from $L^2(M,L;\mu)$:
\[
    \langle s,t\rangle_H
    =
    \int_M h_x(s(x),t(x))\,d\mu(x).
\]
We assume that point evaluations are well-defined on $H$. Equivalently, throughout this finite-dimensional discussion we fix pointwise defined measurable representatives of the elements of $H$, and all identities involving arbitrary $L^2$-sections are understood $\mu$-almost everywhere. This will be the case in the geometric applications, where $H$ is a space of holomorphic sections. For $x\in M$, let
\[
    \operatorname{ev}_x:\left\lbrace\begin{array}{ccc}
    H & \longrightarrow & L_x\\
    s & \longmapsto & s(x),
    \end{array}\right.    
\]
be the evaluation map. Its adjoint $\operatorname{ev}_x^*:L_x\longrightarrow H$ is characterized by
\[
    h_x(s(x),v)
    =
    \langle s,\operatorname{ev}_x^*v\rangle_H,
    \qquad s\in H,\ v\in L_x.
\]
The projection kernel associated with $H$ is the map $B_H(x,y)=\operatorname{ev}_x\operatorname{ev}_y^*\in \operatorname{Hom}(L_y,L_x).$ Equivalently, if $s_1,\ldots,s_N$ is an orthonormal basis of $H$, then
\[
    B_H(x,y)
    =
    \sum_{\ell=1}^N s_\ell(x)\otimes s_\ell(y)^*
    \in L_x\otimes L_y^*
    \simeq \operatorname{Hom}(L_y,L_x),
\]
where $s_\ell(y)^*(v)=h_y(v,s_\ell(y)).$ This kernel is the Schwartz kernel of the orthogonal projection $\Pi_H:L^2(M,L;\mu)\longrightarrow H.$ Indeed, for every $s\in L^2(M,L;\mu)$,
\[
(\Pi_Hs)(x)=\int_M B_H(x,y)s(y)\,d\mu(y)
\]
for $\mu$-almost every $x$; if $s\in H$, the identity holds at every point where the chosen representatives are evaluated. The reproducing identity takes the intrinsic form
\[
    \int_M B_H(x,z)B_H(z,y)\,d\mu(z)
    =
    B_H(x,y),
\]
where the product denotes the composition of operators. On the diagonal, $B_H(x,x)\in \operatorname{End}(L_x)$. Since $L_x$ is one-dimensional, $B_H(x,x)$ is multiplication by a non-negative real number; we denote this scalar by the same symbol. It satisfies
\[
    \int_M B_H(x,x)\,d\mu(x)=N.
\]

The expression $\det(B(x_i,x_j))_{1\leq i,j\leq m}$, which will later characterize the $m$-point correlation function of the Bergman DPP, can be formally written as:
\begin{equation}\label{eq:Bergman_det_section}
{\det}(B(x_i,x_j))=\sum_{\sigma\in\mathfrak{S}_m}\varepsilon(\sigma)\bigotimes\limits_{i=1}^m B(x_i,x_{\sigma(i)})=\sum_{\sigma\in\mathfrak{S}_m}\varepsilon(\sigma)\sum_{i_1,\ldots,i_m=1}^N\bigotimes\limits_{j=1}^m s_{i_j}(x_j)\otimes s_{i_j}(x_{\sigma(j)})^*.\nonumber
\end{equation}
Since it requires adding tensors of different vector spaces, this equation does not make sense algebraically. At best, using isomorphisms between primal and dual fibers by contraction, one could rearrange the terms to define the determinant by
\begin{equation}\label{eq:Bergman_det_function}
{\det}(B(x_i,x_j))=\sum_{\sigma\in\mathfrak{S}_m}\varepsilon(\sigma)\sum_{i_1,\ldots,i_m=1}^N\prod\limits_{j=1}^m h_{x_j}(s_{i_j}(x_j),s_{i_{\sigma^{-1}(j)}}(x_j)).
\end{equation}
This definition is intuitive but not intrinsic, and we will propose another one below.

Let $\mathbf x=(x_1,\ldots,x_m)\in M^m$. Define the finite-dimensional Hermitian vector space
\[
    L_{\mathbf x}
    =
    L_{x_1}\oplus\cdots\oplus L_{x_m}.
\]
The kernel $B_H$ defines an endomorphism $\mathsf B_H(\mathbf x):L_{\mathbf x}\longrightarrow L_{\mathbf x}$ whose $(i,j)$-block is
\[
    B_H(x_i,x_j):L_{x_j}\longrightarrow L_{x_i}.
\]
We define the intrinsic determinant of the line-bundle-valued matrix
$\big(B_H(x_i,x_j)\big)_{1\leq i,j\leq m}$ by
\begin{equation}\label{eq:intrinsic-det}
\det_{\mathbf x}\big(B_H(x_i,x_j)\big)_{1\leq i,j\leq m}:=    \det_{L_{\mathbf x}}\big(\mathsf B_H(\mathbf x)\big).
\end{equation}
This is an ordinary determinant of an endomorphism of the $m$-dimensional complex vector space $L_{\mathbf x}$, and is therefore a well-defined scalar. Since $\mathsf B_H(\mathbf x)=\operatorname{ev}_{\mathbf x}\operatorname{ev}_{\mathbf x}^*,$ where
\[
\operatorname{ev}_{\mathbf x}:\left\lbrace
\begin{array}{lll}
H & \longrightarrow & L_{\mathbf x}\\
s & \longmapsto &  (s(x_1),\ldots,s(x_m)),
\end{array}
\right.
\]
the operator $\mathsf B_H(\mathbf x)$ is non-negative.

\begin{lemma}\label{lem:local-expression-intrinsic-det}
Let $U\subset M$ be an open set over which $L$ admits a non-vanishing local frame $e$. For $x,y\in U$, define the scalar kernel $K_H^e(x,y)$ by
\[
    B_H(x,y)e(y)=K_H^e(x,y)e(x).
\]
Then, for $x_1,\ldots,x_m\in U$,
\[
    \det_{\mathbf x}\big(B_H(x_i,x_j)\big)_{1\leq i,j\leq m}
    =
    \det\big(K_H^e(x_i,x_j)\big)_{1\leq i,j\leq m}.
\]
Moreover, the right-hand side is independent of the choice of the local frame.
\end{lemma}

\begin{proof}
The first equality follows by writing the endomorphism $\mathsf B_H(\mathbf x)$ in the basis $e(x_1),\ldots,e(x_m)$ of $L_{\mathbf x}$. If $e'=ge$ is another local frame, then
\[
    K_H^{e'}(x,y)
    =
    \frac{g(y)}{g(x)}K_H^e(x,y).
\]
Hence
\[
    \big(K_H^{e'}(x_i,x_j)\big)_{i,j}
    =
    D^{-1}\big(K_H^e(x_i,x_j)\big)_{i,j}D,
\]
where $D=\operatorname{diag}(g(x_1),\ldots,g(x_m))$. The determinant is therefore unchanged.
\end{proof}

Lemma~\ref{lem:local-expression-intrinsic-det} is the simplest way to pass between the geometric and the probabilistic languages. Locally, a line-bundle-valued Bergman kernel becomes a scalar kernel, and all usual determinantal formulas apply. However, the intrinsic construction avoids dealing with problematic situations such as the following: take $U,V\subset M$ two disjoint open sets with two separate local frames, then the expression
\[
B_H(x,x)B_H(y,y)-B_H(x,y)B_H(y,x),\quad x\in U,\ y\in V
\]
cannot be easily rewritten using local scalar kernels, which reinforces the choice of taking~\eqref{eq:intrinsic-det} as a definition of the determinant of the kernel.

\begin{proposition}
The intrinsic definition of $\det\limits_{\mathbf x}(B_H(x_i,x_j))$ and its realization~\eqref{eq:Bergman_det_function} are equivalent.
\end{proposition}

\begin{proof}
Let $\mathbf x=(x_1,\ldots,x_m)\in M^m$, and let $s_1,\ldots,s_N$ be an orthonormal basis of $H$. For each $j$, choose a unit vector $\eta_j\in L_{x_j}$. We regard $\eta_1,\ldots,\eta_m$ as an orthonormal basis of the Hermitian vector space $L_{\mathbf x}=L_{x_1}\oplus\cdots\oplus L_{x_m}.$ Write
\[
a_{j\ell}=h_{x_j}(s_\ell(x_j),\eta_j),
\qquad 1\leq j\leq m,\quad 1\leq \ell\leq N.
\]
Since the Hermitian product is linear in the first variable, one has $s_\ell(x_j)=a_{j\ell}\eta_j.$ In the basis $\eta_1,\ldots,\eta_m$, the endomorphism $B_H(\mathbf x)=\operatorname{ev}_{\mathbf x}\operatorname{ev}_{\mathbf x}^*\in \operatorname{End}(L_{\mathbf x})$ has matrix entries $G_{ij}=h_{x_i}\bigl(B_H(x_i,x_j)\eta_j,\eta_i\bigr).$ Using the formula
\[
B_H(x_i,x_j)v=\sum_{\ell=1}^Ns_\ell(x_i)\, h_{x_j}(v,s_\ell(x_j)),\qquad v\in L_{x_j},
\]
we obtain $G_{ij}=\sum_{\ell=1}^Na_{i\ell}\overline{a_{j\ell}}.$ Thus, if $A=(a_{j\ell})_{1\leq j\leq m,\,1\leq \ell\leq N}$, then $G=AA^*.$ Therefore
\[
\det_{\mathbf x}\bigl(B_H(x_i,x_j)\bigr)_{1\leq i,j\leq m}
=
\det(G)=
\sum_{\sigma\in S_m}\varepsilon(\sigma)
\sum_{\ell_1,\ldots,\ell_m=1}^N
\prod_{i=1}^m
a_{i\ell_i}\overline{a_{\sigma(i)\ell_i}}.
\]
Reindexing the conjugated factors by writing $j=\sigma(i)$, this becomes
\[
\det_{\mathbf x}\bigl(B_H(x_i,x_j)\bigr)_{1\leq i,j\leq m}
=\sum_{\sigma\in S_m}\varepsilon(\sigma)\sum_{\ell_1,\ldots,\ell_m=1}^N\prod_{j=1}^ma_{j\ell_j}\overline{a_{j\ell_{\sigma^{-1}(j)}}}.
\]
This is exactly the realization~\eqref{eq:Bergman_det_function}. Since the left-hand side is the intrinsic determinant of the endomorphism of $L_{\mathbf x}$, the expression is independent of all auxiliary choices. The proof is complete.
\end{proof}

\subsection{Projection DPPs associated with sections}

Let $s_1,\ldots,s_N$ be an orthonormal basis of $H$. The analogue of the usual Slater determinant is the section of the external tensor product $L^{\boxtimes N}=\operatorname{pr}_1^*L\otimes\cdots\otimes \operatorname{pr}_N^*L\longrightarrow M^N$ defined by
\[
    \mathfrak S_H(x_1,\ldots,x_N)
    =
    \sum_{\sigma\in\mathfrak S_N}
    \operatorname{sgn}(\sigma)\,
    s_{\sigma(1)}(x_1)\otimes\cdots\otimes s_{\sigma(N)}(x_N).
\]
Thus $\mathfrak S_H(x_1,\ldots,x_N)\in L_{x_1}\otimes\cdots\otimes L_{x_N}.$ Its squared norm is a scalar function on $M^N$, and it is independent of the chosen orthonormal basis of $H$, because a unitary change of basis multiplies $\mathfrak S_H$ by a complex number of modulus one.


We may now define the projection DPP associated with $H$ exactly as in the scalar case.

\begin{definition}\label{def:section-projection-DPP}
Let $H\subset L^2(M,L;\mu)$ be an $N$-dimensional Hilbert space of sections. The projection DPP associated with $H$ is the probability measure on $M^N$ given by
\[
    d\mathbb P_H(x_1,\ldots,x_N)
    =
    \frac{1}{N!}
    \|\mathfrak S_H(x_1,\ldots,x_N)\|^2
    \,d\mu(x_1)\cdots d\mu(x_N).
\]
\end{definition}

The normalization is a direct consequence of the orthonormality of $s_1,\ldots,s_N$, or equivalently of Andreief's identity in local trivializations. The main result of this section is that this projection DPP is genuinely a determinantal point process, which might seem tautological but is not a consequence of existing results on usual projection DPPs because the projection kernel is not scalar. However, the proof will follow from similar arguments.

\begin{theorem}\label{thm:intrinsic-determinantal-property}
Let $X$ be the point process associated with $\mathbb P_H$. Then $X$ is determinantal: for every $1\leq m\leq N$, its $m$-point correlation function is
\[
    \rho_m(x_1,\ldots,x_m)
    =
    \det_{\mathbf x}
    \big(B_H(x_i,x_j)\big)_{1\leq i,j\leq m}.
\]
For $m>N$, one has $\rho_m=0$.
\end{theorem}

\begin{proof}
Let $s_1,\ldots,s_N$ be an orthonormal basis of $H$. For every subset $I=\{i_1<\cdots<i_r\}\subset\{1,\ldots,N\}$, define the partial Slater section
\[
S_I(z_1,\ldots,z_r)
=
\sum_{\sigma\in S_r}
\operatorname{sgn}(\sigma)\,
s_{i_{\sigma(1)}}(z_1)\otimes\cdots\otimes s_{i_{\sigma(r)}}(z_r).
\]
Thus $S_{\{1,\ldots,N\}}=S_H$.

We first record two elementary identities. Let $\mathbf x=(x_1,\ldots,x_m)$. By the Cauchy--Binet formula applied to the evaluation map $\operatorname{ev}_{\mathbf x}:H\longrightarrow
L_{x_1}\oplus\cdots\oplus L_{x_m},$ one has
\begin{equation}\label{eq:32}
\det_{\mathbf x}\bigl(B_H(x_i,x_j)\bigr)_{1\leq i,j\leq m}
=
\sum_{\substack{I\subset\{1,\ldots,N\}\\ |I|=m}}
\|S_I(\mathbf x)\|^2.
\end{equation}
Indeed, after choosing unit vectors in the fibers $L_{x_1},\ldots,L_{x_m}$, the matrix of $B_H(\mathbf x)=\operatorname{ev}_{\mathbf x}\operatorname{ev}_{\mathbf x}^*$ is $AA^*$, where $A$ is the $m\times N$ matrix of the evaluation map. Hence
\[
\det(AA^*)
=
\sum_{|I|=m} |\det A_I|^2,
\]
and $|\det A_I|^2$ is precisely $\|S_I(\mathbf x)\|^2$.

Next, for subsets $I,J\subset\{1,\ldots,N\}$ with $|I|=|J|=r$, orthonormality gives
\begin{equation}\label{eq:33}
\int_{M^r}
\left\langle
S_I(\mathbf z),S_J(\mathbf z)
\right\rangle
\,d\mu(z_1)\cdots d\mu(z_r)
=
r!\,\mathbf 1_{\{I=J\}}.
\end{equation}
This follows by expanding both Slater sections and using
\[
\int_M h_x(s_a(x),s_b(x))\,d\mu(x)=\delta_{ab}.
\]

Now fix $1\leq m\leq N$, and write
\[
\mathbf x=(x_1,\ldots,x_m),
\qquad
\mathbf y=(y_1,\ldots,y_{N-m}).
\]
The Laplace expansion of the full Slater section with respect to the first $m$ variables gives
\[
S_H(\mathbf x,\mathbf y)
=
\sum_{\substack{I\subset\{1,\ldots,N\}\\ |I|=m}}
\varepsilon(I)\,
S_I(\mathbf x)\otimes S_{I^c}(\mathbf y),
\]
where $\varepsilon(I)\in\{\pm1\}$ is the sign of the shuffle putting the indices of $I$ in the first $m$ positions and the indices of $I^c$ in the remaining positions.

Taking the squared norm, integrating over $\mathbf y$, and using~\eqref{eq:33}, we obtain
\[
\begin{aligned}
\int_{M^{N-m}}
\|S_H(\mathbf x,\mathbf y)\|^2
\,d\mu(y_1)\cdots d\mu(y_{N-m})
&=
\sum_{I,J}
\varepsilon(I)\varepsilon(J)
\left\langle S_I(\mathbf x),S_J(\mathbf x)\right\rangle
\\
&\qquad\qquad\times
\int_{M^{N-m}}
\left\langle
S_{I^c}(\mathbf y),S_{J^c}(\mathbf y)
\right\rangle
\,d\mu^{\otimes(N-m)}(\mathbf y)
\\
&=
(N-m)!
\sum_{\substack{I\subset\{1,\ldots,N\}\\ |I|=m}}
\|S_I(\mathbf x)\|^2.
\end{aligned}
\]
Combining this with~\eqref{eq:32}, we get
\begin{equation}\label{eq:35}
\frac{1}{(N-m)!}
\int_{M^{N-m}}
\|S_H(\mathbf x,\mathbf y)\|^2
\,d\mu^{\otimes(N-m)}(\mathbf y)
=
\det_{\mathbf x}\bigl(B_H(x_i,x_j)\bigr)_{1\leq i,j\leq m}.
\end{equation}
Since the density of $P_H$ is $\frac{1}{N!}\|S_H(x_1,\ldots,x_N)\|^2$ with respect to $\mu^{\otimes N}$, the $m$-point correlation function is
\[
\begin{aligned}
\rho_m(x_1,\ldots,x_m)
&=
\frac{N!}{(N-m)!}
\int_{M^{N-m}}
\frac{1}{N!}
\|S_H(x_1,\ldots,x_m,y_1,\ldots,y_{N-m})\|^2
\,d\mu(y_1)\cdots d\mu(y_{N-m})
\\
&=
\frac{1}{(N-m)!}
\int_{M^{N-m}}
\|S_H(\mathbf x,\mathbf y)\|^2
\,d\mu^{\otimes(N-m)}(\mathbf y).
\end{aligned}
\]
By~\eqref{eq:35}, this equals
\[
\rho_m(x_1,\ldots,x_m)
=
\det_{\mathbf x}\bigl(B_H(x_i,x_j)\bigr)_{1\leq i,j\leq m}.
\]
This proves the claimed determinantal formula for $1\leq m\leq N$.

Finally, since the process has exactly $N$ particles, its factorial moment measures of order $m>N$ vanish. Hence $\rho_m=0$ for $m>N$.
\end{proof}

It is important to emphasize that, although the kernel $B_H$ is not a scalar function on $M\times M$, the correlation functions are genuine scalar functions on $M^m$. It follows that we have a determinantal point process in the usual probabilistic sense. The point of the theorem is not merely that one may trivialize the line bundle locally and apply the scalar theory: local scalar kernels depend on the chosen frame, and the corresponding matrices are conjugated by diagonal transition matrices. The intrinsic determinant is the invariant scalar left by this conjugation. Thus the DPP itself lives on $M$, not on a choice of local gauge. This is the probabilistic analogue of the usual gauge invariance of Bergman kernels and of fermionic correlation functions.

\begin{remark}
If $L=M\times\mathbb C$ is the trivial line bundle with its standard Hermitian metric, then $L_x=\mathbb C$ for every $x$, and $\operatorname{Hom}(L_y,L_x)\simeq \mathbb C.$ The intrinsic determinant reduces to the usual determinant:
\[
    \det_{\mathbf x}
    \big(B_H(x_i,x_j)\big)_{1\leq i,j\leq m}
    =
    \det\big(K_H(x_i,x_j)\big)_{1\leq i,j\leq m},
\]
and we recover the usual construction of finite-rank scalar projection DPPs.
\end{remark}

To conclude this section, let us describe the local form of the scalar kernel introduced in Lemma~\ref{lem:local-expression-intrinsic-det}. Let $e$ be a local frame of $L$ on $U$, and write $\|e(x)\|_h^2=e^{-\varphi(x)}.$ If $s_\ell=f_\ell e$, then
\[
    B_H(x,y)e(y)
    =
    \left(
        \sum_{\ell=1}^N
        f_\ell(x)\overline{f_\ell(y)}e^{-\varphi(y)}
    \right)e(x).
\]
Thus
\[
    K_H^e(x,y)
    =
    \sum_{\ell=1}^N
    f_\ell(x)\overline{f_\ell(y)}e^{-\varphi(y)}.
\]
If instead one uses the local unit frame $\widetilde e(x)=e^{\varphi(x)/2}e(x),$ then the scalar kernel becomes
\[
    \widetilde K_H(x,y)
    =
    e^{-\varphi(x)/2}
    \left(
        \sum_{\ell=1}^N f_\ell(x)\overline{f_\ell(y)}
    \right)
    e^{-\varphi(y)/2}.
\]
This is the form in which the kernel most closely resembles the usual Christoffel--Darboux kernel. In particular, in local unitary frames the line-bundle-valued formalism becomes indistinguishable from the scalar projection-DPP formalism of the previous section.



\section{Transfer principles from kernels to DPPs}
\label{sec:black-box-transfer}

We now isolate a simple principle that will be used throughout the rest of the paper: DPPs are controlled by their kernels. Consequently, many probabilistic statements about projection DPPs follow from analytic estimates on their projection kernels. The results presented here are not meant as new structural theorems on determinantal point processes; in fact, most of them are standard consequences of the finite-rank projection-DPP formalism. Their role here is different: we isolate the exact finite-dimensional mechanisms by which analytic information on a sequence of projection kernels is converted into probabilistic information on the associated point processes.


Let $E$ be a compact metric space. For every $k\geq 1$, let $\nu_k$ be a finite positive Borel measure on $E$, and let $H_k\subset L^2(E,\nu_k)$ be a finite-dimensional subspace of dimension $N_k$. We assume throughout that $N_k\longrightarrow\infty.$ Let $K_k$ be the kernel of the orthogonal projection $\Pi_k:L^2(E,\nu_k)\longrightarrow H_k.$ The associated projection DPP will be denoted by $X_k=\{X_{1,k},\ldots,X_{N_k,k}\}.$ Its $m$-point correlation functions with respect to $\nu_k$ are
\[
    \rho_{m,k}(x_1,\ldots,x_m)
    =
    \det\big(K_k(x_i,x_j)\big)_{1\leq i,j\leq m},
    \qquad 1\leq m\leq N_k.
\]

The empirical measure of $X_k$ is $\widehat\mu_k=\frac{1}{N_k}\sum_{x\in X_k}\delta_x.$ We shall also use the normalized one-point intensity measure $\beta_k=\frac{1}{N_k}K_k(x,x)\,d\nu_k(x).$


\subsection{Macroscopic transfer: diagonal and off-diagonal estimates}

The first transfer principle says that convergence of the normalized diagonal of the kernel implies convergence of the empirical measure.

\begin{theorem}\label{thm:diagonal-transfer}
Assume that the normalized one-point intensity measures satisfy
\[
    \beta_k
    =
    \frac{1}{N_k}K_k(x,x)\,d\nu_k(x)
    \Longrightarrow
    \nu
\]
weakly, for some probability measure $\nu$ on $E$. Then $\widehat\mu_k\Longrightarrow\nu$ in probability, with respect to the weak topology on probability measures on
$E$.
\end{theorem}

\begin{proof}
Since $E$ is compact, it is enough to prove that, for every continuous function $f\in C^0(E)$,
\[
    \widehat\mu_k(f)
    =
    \frac{1}{N_k}\sum_{x\in X_k}f(x)
\]
converges in probability to $\nu(f)$.

By the one-point correlation formula,
\[
    \mathbb E[\widehat\mu_k(f)]
    =
    \frac{1}{N_k}\int_E f(x)K_k(x,x)\,d\nu_k(x)
    =
    \int_E f\,d\beta_k.
\]
By assumption, this converges to $\int_E f\,d\nu$.

It remains to control the variance. Using the variance identity~\eqref{eq:variance-dpp},
\[
    \operatorname{Var}\big(\widehat\mu_k(f)\big)
    =
    \frac{1}{2N_k^2}
    \int_{E^2}
        |f(x)-f(y)|^2|K_k(x,y)|^2
        \,d\nu_k(x)d\nu_k(y).
\]
Since $|f(x)-f(y)|^2\leq 4\|f\|_\infty^2$ and $K_k$ is a projection kernel,
\[
    \int_{E^2}|K_k(x,y)|^2\,d\nu_k(x)d\nu_k(y)
    =
    \operatorname{Tr}(\Pi_k^2)
    =
    \operatorname{Tr}(\Pi_k)
    =
    N_k.
\]
Thus $\operatorname{Var}\big(\widehat\mu_k(f)\big)\leq\frac{2\|f\|_\infty^2}{N_k}\longrightarrow 0$, and the convergence in probability follows.
\end{proof}

For later use, we record the exact variance identity in a form that is often more informative than the crude $O(N_k^{-1})$ bound, which uses also~\eqref{eq:variance-dpp}.

\begin{proposition}\label{prop:variance-transfer}
Let $f$ be a bounded real-valued measurable function on $E$. Then
\[
    \operatorname{Var}\big(\widehat\mu_k(f)\big)
    =
    \frac{1}{2N_k^2}
    \int_{E^2}
        |f(x)-f(y)|^2
        |K_k(x,y)|^2
        \,d\nu_k(x)d\nu_k(y).
\]
In particular, if the mass of $|K_k(x,y)|^2d\nu_k(x)d\nu_k(y)$ is concentrated near the diagonal, then the variance of smooth linear statistics is controlled by the modulus of continuity of $f$.
\end{proposition}

\begin{corollary}\label{cor:variance-off-diagonal}
Let $d$ be a metric on $E$. Assume that there exists a sequence
$\varepsilon_k\to0$ such that
\[
    \frac{1}{N_k}
    \int_{\{d(x,y)>\varepsilon_k\}}
        |K_k(x,y)|^2
        \,d\nu_k(x)d\nu_k(y)
    \longrightarrow 0.
\]
Then, for every continuous $f\in C^0(E)$,
\[
    \operatorname{Var}\big(\widehat\mu_k(f)\big)
    \leq
    \frac12\omega_f(\varepsilon_k)^2\frac{1}{N_k}
    +
    \frac{2\|f\|_\infty^2}{N_k^2}
    \int_{\{d(x,y)>\varepsilon_k\}}
        |K_k(x,y)|^2
        \,d\nu_k(x)d\nu_k(y),
\]
where
\[
    \omega_f(\varepsilon)
    =
    \sup\{|f(x)-f(y)|:d(x,y)\leq\varepsilon\}.
\]
\end{corollary}

\begin{proof}
Split the variance integral into the regions $d(x,y)\leq\varepsilon_k$ and $d(x,y)>\varepsilon_k$. On the first region,
\[
    |f(x)-f(y)|^2\leq \omega_f(\varepsilon_k)^2,
\]
and
\[
    \int_{E^2}|K_k(x,y)|^2\,d\nu_k(x)d\nu_k(y)=N_k.
\]
On the second region, use
\[
    |f(x)-f(y)|^2\leq 4\|f\|_\infty^2.
\]
This gives the claimed bound.
\end{proof}


\subsection{Microscopic transfer: local kernels and Palm kernels}\label{sec:Palm-transfer}

The next transfer principle explains why local asymptotics of kernels imply local limits of DPP correlations. This is the basic mechanism behind Ginibre or Bargmann--Fock universality.

Let $p\in E$ be a point around which $E$ is locally modeled by an open subset of $\mathbb R^d$. In the geometric applications, $d=2n$ and the local model will be $\mathbb C^n$. Let
\[
    \chi:U\longrightarrow V\subset\mathbb R^d
\]
be a coordinate chart with $\chi(p)=0$. Let $r_k\to0$ be a microscopic scale, and define
\[
    x_k(u)=\chi^{-1}(r_k u)
\]
whenever $u$ belongs to a fixed compact subset of $\mathbb R^d$.

We formulate the statement after writing the kernel with respect to Lebesgue measure in the chosen coordinates. Thus, if $\nu_k$ has density $a_k(z)\,dz$ in $U$, and if $K_k^U(z,w)$ denotes the scalar kernel with respect to $\nu_k$, we set
\[
    \mathcal K_k^U(z,w)
    =
    a_k(z)^{1/2}K_k^U(z,w)a_k(w)^{1/2}.
\]
The local correlation functions with respect to Lebesgue measure are then determinants of $\mathcal K_k^U$.

\begin{theorem}\label{thm:local-kernel-transfer}
Let $r_k\to0$ and $A_k\to\infty$. Assume that, locally uniformly for
$u,v\in\mathbb R^d$,
\[
    A_k^{-1}
    \mathcal K_k^U(r_k u,r_k v)
    \longrightarrow
    K_\infty(u,v),
\]
where $K_\infty$ is a continuous Hermitian kernel on $\mathbb R^d$. Then, for every fixed $m\geq1$,
\[
    A_k^{-m}
    \rho_{m,k}^{\mathrm{loc}}(r_k u_1,\ldots,r_k u_m)
    \longrightarrow
    \det\big(K_\infty(u_i,u_j)\big)_{1\leq i,j\leq m},
\]
locally uniformly for $(u_1,\ldots,u_m)\in(\mathbb R^d)^m$.

Moreover, if $r_k^d A_k\longrightarrow c\in(0,\infty),$ then the $m$-point correlation density of the rescaled process $X_k^{(p)}=\sum_{x\in X_k\cap U}\delta_{\chi(x)/r_k}$ with respect to Lebesgue measure on $(\mathbb R^d)^m$ converges locally uniformly to $\det\big(cK_\infty(u_i,u_j)\big)_{1\leq i,j\leq m}.$
\end{theorem}

\begin{proof}
By the determinantal formula,
\[
    \rho_{m,k}^{\mathrm{loc}}(z_1,\ldots,z_m)
    =
    \det\big(\mathcal K_k^U(z_i,z_j)\big)_{1\leq i,j\leq m}.
\]
Taking $z_i=r_ku_i$, we get
\[
    A_k^{-m}\rho_{m,k}^{\mathrm{loc}}(r_ku_1,\ldots,r_ku_m)
    =
    \det\big(A_k^{-1}\mathcal K_k^U(r_ku_i,r_ku_j)\big)_{1\leq i,j\leq m}.
\]
The determinant is a polynomial in the entries of the matrix, hence locally uniform convergence of the entries implies locally uniform convergence of the determinants.

For the rescaled process, the Jacobian of the change of variables gives
\[
    \rho_{m,k}^{(p)}(u_1,\ldots,u_m)
    =
    r_k^{dm}\rho_{m,k}^{\mathrm{loc}}(r_ku_1,\ldots,r_ku_m).
\]
Combining this identity with the first part gives the claimed limit when $r_k^dA_k\to c$.
\end{proof}

The preceding theorem is a correlation-level statement. For convergence of the point processes themselves it is convenient to state the standard trace-norm condition on the rescaled kernels. Define
\[
    \mathcal K_k^{(p)}(u,v)
    =
    r_k^d\mathcal K_k^U(r_ku,r_kv).
\]

\begin{proposition}\label{prop:local-process-transfer}
Assume that $r_k^dA_k\to c\in(0,\infty)$, that the assumptions of Theorem~\ref{thm:local-kernel-transfer} hold, and that the kernel $cK_\infty$ defines a locally trace-class DPP on $\mathbb R^d$. Assume moreover that, for every compact $C\subset\mathbb R^d$, the integral operators on $L^2(C)$ with kernels
\[
    \mathbf 1_C(u)\mathcal K_k^{(p)}(u,v)\mathbf 1_C(v)
\]
converge in trace norm to the operator with kernel
\[
    \mathbf 1_C(u)cK_\infty(u,v)\mathbf 1_C(v).
\]
Then $X_k^{(p)}$ converges in distribution, for the vague topology on locally finite configurations in $\mathbb R^d$, to the DPP with kernel $cK_\infty$.
\end{proposition}

\begin{proof}
Let $\phi$ be a compactly supported measurable function with $0\leq\phi\leq1$. For a finite-rank DPP,
\[
    \mathbb E\left[\prod_{x\in X_k^{(p)}}(1-\phi(x))\right]
    =
    \det\bigl(I-\sqrt{\phi}\,\mathcal K_k^{(p)}\,\sqrt{\phi}\bigr),
\]
where the determinant is the Fredholm determinant on a compact set containing the support of $\phi$. The trace-norm convergence of the restricted kernels implies convergence of these Fredholm determinants to $\det\bigl(I-\sqrt{\phi}\,cK_\infty\,\sqrt{\phi}\bigr),$ which is the corresponding multiplicative functional of the limiting DPP. These multiplicative functionals determine the law of a simple locally finite point process, and hence give convergence in distribution.
\end{proof}

We can also go beyond the mere convergence of correlation functions: assume that, for some sequence \(\varepsilon_k\to0\), one has locally uniformly on \(\mathbb R^d\times\mathbb R^d\)
\[
A_k^{-1}\mathcal K_k(r_ku,r_kv) =\sum_{j=0}^N \varepsilon_k^j K_j(u,v)+O(\varepsilon_k^{N+1})
\]
for every \(N\ge0\). Then, since the determinant is polynomial in matrix entries, for every fixed \(m\ge1\),
\begin{equation}\label{eq:expansion-correl}
A_k^{-m}\rho_{m,k}^{\mathrm{loc}}(r_ku_1,\ldots,r_ku_m)=\sum_{j=0}^N \varepsilon_k^j R_{m,j}(u_1,\ldots,u_m)+O(\varepsilon_k^{N+1}),
\end{equation}
locally uniformly on \((\mathbb R^d)^m\), where
\[
R_{m,j}=[t^j]\det\left(K_0(u_a,u_b)+\sum_{q\ge1}t^qK_q(u_a,u_b)\right)_{1\le a,b\le m}.
\]
In particular,
\[
R_{m,0}(u_1,\ldots,u_m)=\det(K_0(u_a,u_b))_{a,b=1}^m.
\]


Now, recall that projection DPPs are stable under Palm conditioning, and their Palm kernels are Schur complements. Therefore, local kernel asymptotics automatically imply local Palm asymptotics, provided the finite matrices at the conditioned points remain invertible.

We use the notation of the previous subsection. Let $q_1,\ldots,q_\ell\in\mathbb R^d$ be distinct points and set
\[
    p_{a,k}=x_k(q_a)=\chi^{-1}(r_kq_a).
\]
Assume that the Gram matrix
\[
    G_\infty(\mathbf q)
    =
    \big(K_\infty(q_a,q_b)\big)_{1\leq a,b\leq \ell}
\]
is invertible.

\begin{proposition}\label{prop:palm-transfer}
Under the assumptions of Theorem~\ref{thm:local-kernel-transfer}, the reduced Palm kernels satisfy
\[
    A_k^{-1}
    \mathcal K_k^{!,\mathbf p_k}(r_ku,r_kv)
    \longrightarrow
    K_\infty^{!,\mathbf q}(u,v)
\]
locally uniformly, where
\[
    K_\infty^{!,\mathbf q}(u,v)
    =
    K_\infty(u,v)
    -
    \sum_{a,b=1}^{\ell}
    K_\infty(u,q_a)
    \big(G_\infty(\mathbf q)^{-1}\big)_{ab}
    K_\infty(q_b,v).
\]
If in addition $r_k^dA_k\to c\in(0,\infty)$, then the rescaled reduced Palm kernels with respect to Lebesgue measure satisfy
\[
    r_k^d\mathcal K_k^{!,\mathbf p_k}(r_ku,r_kv)
    \longrightarrow
    cK_\infty^{!,\mathbf q}(u,v).
\]
Consequently, the local correlation functions of the reduced Palm processes converge to those of the reduced Palm DPP with kernel $cK_\infty^{!,\mathbf q}$.
\end{proposition}

\begin{proof}
The finite-$k$ reduced Palm kernel is the Schur complement
\[
    \mathcal K_k^{!,\mathbf p_k}(z,w)
    =
    \mathcal K_k^U(z,w)
    -
    \sum_{a,b=1}^{\ell}
    \mathcal K_k^U(z,r_kq_a)
    \big(G_k(\mathbf q)^{-1}\big)_{ab}
    \mathcal K_k^U(r_kq_b,w),
\]
where $G_k(\mathbf q)=\big(\mathcal K_k^U(r_kq_a,r_kq_b)\big)_{1\leq a,b\leq \ell}.$ By the local kernel convergence,
\[
    A_k^{-1}G_k(\mathbf q)\longrightarrow G_\infty(\mathbf q).
\]
Since $G_\infty(\mathbf q)$ is invertible,
\[
    A_kG_k(\mathbf q)^{-1}\longrightarrow G_\infty(\mathbf q)^{-1}.
\]
Substituting these limits in the Schur complement formula gives the announced convergence. The convergence of Palm correlation functions then follows by taking determinants and by using the same Jacobian factor as in Theorem~\ref{thm:local-kernel-transfer}.
\end{proof}

\subsection{Laplace transform transfer: Toeplitz traces and large deviations}

Another important transfer principle concerns multiplicative functionals. It will later explain why changes of one-body weights in Bergman ensembles are encoded by determinants of compressed multiplication operators.

Let $f:E\to\mathbb R$ be bounded. For the projection DPP with kernel $K_k$, define $S_{f,k}=\sum_{x\in X_k}f(x).$ Let $T_{g,k}=\Pi_k M_g\Pi_k$ be the compression to $H_k$ of multiplication by a bounded function $g$. The following is a direct consequence of~\eqref{eq:exp-lin-stat}.

\begin{proposition}\label{prop:multiplicative-transfer}
For every bounded measurable function $f:E\to\mathbb R$,
\begin{equation}
\mathbb E\left[e^{S_{f,k}}\right]=\det_{H_k}\big(I+T_{e^f-1,k}\big).
\end{equation}
\end{proposition}


\begin{remark}
If one has asymptotics for traces of products of compressed multiplication operators, then Proposition~\ref{prop:multiplicative-transfer} gives asymptotics for the Laplace transforms of linear statistics. In the Bergman setting, the operators $T_{g,k}$ are Toeplitz operators. Thus Toeplitz trace asymptotics become probabilistic limit theorems for Bergman DPPs.
\end{remark}

For completeness, we record the finite-$k$ cumulant formula. It will be used later only as a guide, but it explains why Toeplitz calculus is naturally adapted to linear statistics of Bergman DPPs. If $f:E\to\mathbb R$ is a bounded function, set $\kappa_{\ell,k}(f)$ for the $\ell$-th cumulant of $S_{f,k}$.

\begin{proposition}\label{prop:cumulant-trace-formula}
For every $\ell\geq1$,
\[
    \kappa_{\ell,k}(f)
    =
    \sum_{r=1}^{\ell}
    \frac{(-1)^{r-1}}{r}
    \sum_{\substack{m_1+\cdots+m_r=\ell\\ m_i\geq1}}
    \frac{\ell!}{m_1!\cdots m_r!}
    \operatorname{Tr}_{H_k}
    \left(
        T_{f^{m_1},k}\cdots T_{f^{m_r},k}
    \right).
\]
\end{proposition}

\begin{proof}
By Proposition~\ref{prop:multiplicative-transfer},
\[
    \log\mathbb E[e^{tS_{f,k}}]
    =
    \operatorname{Tr}
    \log\big(I+T_{e^{tf}-1,k}\big).
\]
Using the power-series expansions $\log(I+A)=\sum_{r\geq1}\frac{(-1)^{r-1}}{r}A^r$ and $e^{tf}-1=\sum_{m\geq1}\frac{t^m}{m!}f^m,$ we identify the coefficient of $t^\ell/\ell!$. Since the space $H_k$ is finite-dimensional, the expansion is purely algebraic in a neighborhood of $t=0$.
\end{proof}

\begin{remark}
This formula separates probability from analysis. The probabilistic identity is exact and finite-dimensional. Any limiting statement about the cumulants is then reduced to asymptotic estimates for traces of products of the operators $T_{g,k}$. For Bergman kernels, these are precisely Berezin--Toeplitz operators.
\end{remark}

The multiplicative functional identity also gives a convenient abstract route to large deviations for empirical measures. The point is that the logarithmic moment-generating function of the empirical measure can be written as a determinant of a compressed multiplication operator.

Let $a_k\to\infty$ be a sequence of positive numbers. For $f\in C(E,\mathbb R)$, define
\[
\Lambda_k(f)
=
\frac{1}{a_k}
\log
\mathbb E
\left[
\exp\left(a_k\widehat\mu_k(f)\right)
\right],
\qquad
\widehat\mu_k(f)=\int_E f\,d\widehat\mu_k
=
\frac{1}{N_k}\sum_{x\in X_k}f(x).
\]
Since $a_k\widehat\mu_k(f)=\frac{a_k}{N_k}S_{f,k},$ Proposition~\ref{prop:multiplicative-transfer} gives
\[
\Lambda_k(f)
=
\frac{1}{a_k}
\log
\det_{H_k}
\left(
I+
T_{\exp((a_k/N_k)f)-1,k}
\right).
\]
Thus large deviations for the empirical measures are reduced to asymptotics of determinants of compressed multiplication operators.

\begin{theorem}\label{thm:ldp-transfer}
Let $E$ be compact, and let $X_k$ be the projection DPP associated with $H_k\subset L^2(E,\nu_k)$. Let $a_k\to\infty$. Assume that, for every $f\in C(E,\mathbb R)$, the limit
\begin{equation}\label{eq:46}
\Lambda(f)
=
\lim_{k\to\infty}
\frac{1}{a_k}
\log
\det_{H_k}
\left(
I+
T_{\exp((a_k/N_k)f)-1,k}
\right)
\end{equation}
exists as a finite real number. Assume moreover that $\Lambda$ is Gateaux differentiable on $C(E,\mathbb R)$. Then the laws of the empirical measures $\widehat\mu_k$ on $\mathcal P(E)$, endowed with the weak topology, satisfy a large-deviation principle with speed $a_k$ and good rate function
\[
I(\nu)
=
\sup_{f\in C(E,\mathbb R)}
\left\{
\int_E f\,d\nu-\Lambda(f)
\right\}.
\]
\end{theorem}

\begin{proof}
The space $\mathcal P(E)$ is compact for the weak topology, since $E$ is compact; hence the sequence of laws of $\widehat\mu_k$ is exponentially tight.

Fix $f_1,\ldots,f_m\in C(E,\mathbb R)$. For $t=(t_1,\ldots,t_m)\in\mathbb R^m$, set $f_t=\sum_{j=1}^m t_jf_j$. By Proposition~\ref{prop:multiplicative-transfer},
\[
\frac{1}{a_k}
\log
\mathbb E
\left[
\exp\left(
a_k\sum_{j=1}^m t_j\widehat\mu_k(f_j)
\right)
\right]
=
\Lambda_k(f_t).
\]
By assumption this converges to $\Lambda(f_t)$. Since $\Lambda$ is Gateaux differentiable on $C(E,\mathbb R)$, the function $t\longmapsto \Lambda(f_t)$ is differentiable on $\mathbb R^m$. The finite-dimensional G\"artner--Ellis theorem \cite[Thm.~2.3.6]{DZ} therefore gives an LDP for $\bigl(\widehat\mu_k(f_1),\ldots,\widehat\mu_k(f_m)\bigr)$ with speed $a_k$ and rate function equal to the Legendre transform of $t\mapsto\Lambda(f_t)$.

Applying the Dawson--G\"artner projective-limit theorem \cite[Thm.~4.6.1]{DZ} to the family of maps
\[
\mathcal P(E)\longrightarrow\mathbb R^m,\qquad
\nu\longmapsto \left(\int_E f_1\,d\nu,\ldots,\int_E f_m\,d\nu\right),
\]
and using that these maps generate the weak topology on $\mathcal P(E)$, yields an LDP on $\mathcal P(E)$. The resulting rate function is
\[
I(\nu)
=
\sup_{\substack{m\geq1\\ f_1,\ldots,f_m\in C(E)\\ t\in\mathbb R^m}}
\left\{
\sum_{j=1}^m t_j\int_E f_j\,d\nu
-
\Lambda\left(\sum_{j=1}^m t_jf_j\right)
\right\}.
\]
Since $\sum_jt_jf_j$ ranges over $C(E,\mathbb R)$, this is exactly
\[
I(\nu)
=
\sup_{f\in C(E,\mathbb R)}
\left\{
\int_E f\,d\nu-\Lambda(f)
\right\}.
\]
Compactness of $\mathcal P(E)$ makes the rate function good.
\end{proof}

\begin{corollary}\label{cor:dev-transfer}
Let $V\in C(E,\mathbb R)$. Define a tilted probability measure by
\[
d\mathbb P_k^V
=
\frac{
\exp\bigl(-a_k\widehat\mu_k(V)\bigr)
}{
\mathbb E\left[\exp\bigl(-a_k\widehat\mu_k(V)\bigr)\right]
}
\,d\mathbb P_k.
\]
Assume that the hypotheses of the preceding theorem hold. Then the empirical measures under $\mathbb P_k^V$ satisfy a large-deviation principle with speed $a_k$ and good rate function
\begin{equation}\label{eq:48}
I_V(\nu)
=
I(\nu)+\int_E V\,d\nu
-
\inf_{\eta\in\mathcal P(E)}
\left(
I(\eta)+\int_E V\,d\eta
\right).
\end{equation}

Equivalently, their limiting log-Laplace functional is
\[
\Lambda^V(f)
=
\Lambda(f-V)-\Lambda(-V),
\]
and the rate function is the Legendre--Fenchel transform of $\Lambda^V$.
\end{corollary}

\begin{proof}
For $f\in C(E,\mathbb R)$,
\[
\frac{1}{a_k}
\log
\mathbb E_{\mathbb P_k^V}
\left[
\exp\left(a_k\widehat\mu_k(f)\right)
\right]
=
\Lambda_k(f-V)-\Lambda_k(-V).
\]
Taking $k\to\infty$ gives $\Lambda^V(f)=\Lambda(f-V)-\Lambda(-V)$. The large-deviation principle follows again from the preceding theorem. The expression~\eqref{eq:48} is the usual tilted form of the rate function, obtained by rewriting the Legendre--Fenchel transform of $\Lambda^V$.
\end{proof}

\section{Bergman ensembles}
\label{sec:Bergman-DPP}

We now specialize the preceding construction to the complex geometric setting: the finite-dimensional Hilbert space $H^0(M,L^k)$ is taken as a space of holomorphic sections of a high tensor power of a positive line bundle, and its projection kernel is the Bergman kernel. There are two sources of weights: the Hermitian metric $h^k$ on $L^k$, whose curvature determines the semiclassical geometry, and the auxiliary measure $d\nu_\psi=e^{-\psi}d\mu$, which determines the $L^2$-orthogonal projection. For fixed $k$ both enter only through the Hilbert structure on $H^0(M,L^k)$. In the limit $k\to\infty$, however, the metric $h^k$ is the dominant geometric object, while a smooth auxiliary density affects only lower-order coefficients. From the probabilistic point of view, no new principle is involved: a Bergman ensemble is a finite-rank projection DPP. The geometric content lies in the special origin of the Hilbert spaces and in the asymptotic information available for their projection kernels.

\subsection{Hermitian line bundle and Bergman kernel}

Let $M$ be a compact complex manifold of complex dimension $n$, and let $L\longrightarrow M$ be a holomorphic line bundle. We denote by $h$ a smooth Hermitian metric on $L$. If $e$ is a non-vanishing local holomorphic frame of $L$ on an open set $U\subset M$, we write $\|e(x)\|_h^2=e^{-\varphi(x)}.$ The function $\varphi$ is called the local weight of $h$ in the frame $e$. With the convention
\[
    d^c=\frac{1}{4\pi i}(\partial-\overline{\partial}),
    \qquad
    dd^c=\frac{i}{2\pi}\partial\overline{\partial},
\]
the first Chern form of $(L,h)$ is locally $\omega_h=c_1(L,h)=dd^c\varphi.$ This is the normalization used below in the asymptotic Riemann--Roch formula and in the leading term of the Bergman kernel expansion. We shall say that $L$ is positive if $h$ may be chosen so that $\omega_h$ is a positive $(1,1)$-form. In that case $\omega_h$ is a K\"ahler form. The positivity assumption is not needed for the finite-dimensional construction of the DPP. It becomes important later, when one studies the asymptotic regime $k\to\infty$, because it is precisely the hypothesis under which Bergman kernels admit universal asymptotic expansions.

Let $\mu$ be a positive Borel measure on $M$. In most geometric applications, $\mu$ will be the measure associated with a smooth positive volume form, but the finite-dimensional construction only requires enough integrability for the inner products below to be well-defined.

For every integer $k\geq 1$, let $L^k=L^{\otimes k}$
and denote by $H_k=H^0(M,L^k)$ the finite-dimensional vector space of holomorphic sections of $L^k$.
We write $N_k=\dim H^0(M,L^k).$ Let $\psi:M\to\mathbb R$ be a bounded measurable function, or more generally a function such that the following integrals are finite. We define the weighted measure
\[
    d\nu_\psi(x)=e^{-\psi(x)}\,d\mu(x).
\]
The metric $h$ on $L$ induces a metric $h^k$ on $L^k$, and we endow $H_k$ with the Hermitian inner product
\[
    \langle s,t\rangle_{k,\psi}
    =
    \int_M h_x^k(s(x),t(x))\,e^{-\psi(x)}\,d\mu(x)
    =
    \int_M h_x^k(s(x),t(x))\,d\nu_\psi(x).
\]
When it is useful to display the local weight of $h$, we shall also write $\langle s,t\rangle_{k\varphi+\psi,\mu}.$ Indeed, if $s=f e^{\otimes k}$ and $t=g e^{\otimes k}$ on a local
trivialization $U$, then locally
\[
    h_x^k(s(x),t(x))e^{-\psi(x)}
    =
    f(x)\overline{g(x)}e^{-k\varphi(x)-\psi(x)}.
\]

\begin{remark}[Weights versus reference measures]
The choice of $\psi$ can be viewed in two equivalent ways. Either we keep the metric $h^k$ fixed and change the reference measure from $\mu$ to $\nu_\psi=e^{-\psi}\mu$, or we keep $\mu$ fixed and tensor $L^k$ with the trivial line bundle endowed with the metric of weight $\psi$. These two points of view lead to the same Hilbert space and the same point process.
\end{remark}

The Bergman projection is the orthogonal projection
\[
    \Pi_{k,\psi}:L^2(M,L^k;\nu_\psi)
    \longrightarrow
    H^0(M,L^k).
\]
Its Schwartz kernel with respect to the measure $\nu_\psi$ is denoted by $B_{k,\psi}(x,y)\in \operatorname{Hom}(L_y^k,L_x^k).$ As a line-bundle-valued projection kernel, it satisfies all properties described in Section~\ref{sec:line-bundle-kernels}.


\begin{definition}
The scalar function $x\longmapsto B_{k,\psi}(x,x)$ is called the \emph{Bergman density} of $H^0(M,L^k)$ with respect to $(h^k,\nu_\psi)$. The probability measure
\[
    \beta_{k,\psi}
    =
    \frac{1}{N_k}B_{k,\psi}(x,x)\,d\nu_\psi(x)
\]
will be called the \emph{normalized Bergman measure}.
\end{definition}

Let us show how this translates in local coordinates: let $e$ be a local holomorphic frame of $L$ on $U\subset M$, and write $\|e(x)\|_h^2=e^{-\varphi(x)}.$ If $s_\ell^{(k)}=f_\ell^{(k)}e^{\otimes k}$ on $U$, then the orthonormality condition reads
\[
    \int_M
    h_x^k(s_i^{(k)}(x),s_j^{(k)}(x))e^{-\psi(x)}\,d\mu(x)
    =
    \delta_{ij},
\]
and locally it takes the form
\[
    \int_U
    f_i^{(k)}(x)\overline{f_j^{(k)}(x)}
    e^{-k\varphi(x)-\psi(x)}
    \,d\mu(x)
    +\text{contribution outside }U
    =
    \delta_{ij}.
\]

In the frame $e^{\otimes k}$, the Bergman kernel is represented by a scalar kernel $K_{k,\psi}^e(x,y)$ defined by
\[
    B_{k,\psi}(x,y)e(y)^{\otimes k}
    =
    K_{k,\psi}^e(x,y)e(x)^{\otimes k}.
\]
Using the expression above, one gets
\[
    K_{k,\psi}^e(x,y)
    =
    \sum_{\ell=1}^{N_k}
    f_\ell^{(k)}(x)\overline{f_\ell^{(k)}(y)}e^{-k\varphi(y)}.
\]
This formula is asymmetric because the frame $e$ is not unitary. If instead one passes formally to the local unit frame $\widetilde e(x)=e^{\varphi(x)/2}e(x),$ then the kernel is represented by
\[
    \widetilde K_{k,\psi}(x,y)
    =
    e^{-k\varphi(x)/2}
    \left(
        \sum_{\ell=1}^{N_k}
        f_\ell^{(k)}(x)\overline{f_\ell^{(k)}(y)}
    \right)
    e^{-k\varphi(y)/2}.
\]

If one wants to write the DPP correlations with respect to the original measure $\mu$, rather than with respect to $\nu_\psi$, the convenient local scalar kernel is
\[
    \widehat K_{k,\psi}(x,y)
    =
    e^{-(k\varphi(x)+\psi(x))/2}
    \left(
        \sum_{\ell=1}^{N_k}
        f_\ell^{(k)}(x)\overline{f_\ell^{(k)}(y)}
    \right)
    e^{-(k\varphi(y)+\psi(y))/2}.
\]
This is the ordinary projection kernel associated with the local orthonormal functions
\[
    f_\ell^{(k)}(x)e^{-(k\varphi(x)+\psi(x))/2}.
\]
Thus the same process may be described either by the bundle-valued kernel $B_{k,\psi}$ with respect to $\nu_\psi$, or by the local scalar kernel $\widehat K_{k,\psi}$ with respect to $\mu$.

\subsection{The determinantal point process}

Let $\{s_1^{(k)},\ldots,s_{N_k}^{(k)}\}$ be an orthonormal basis of $H_k$ for $\langle\cdot,\cdot\rangle_{k,\psi}$. Define the Slater determinant section
\begin{equation}\label{eq:Slater}
\mathfrak S_{k,\psi}(x_1,\ldots,x_{N_k})
    =
    \sum_{\sigma\in\mathfrak S_{N_k}}
    \operatorname{sgn}(\sigma)\,
    s_{\sigma(1)}^{(k)}(x_1)\otimes\cdots\otimes
    s_{\sigma(N_k)}^{(k)}(x_{N_k}).
\end{equation}
It is a section of $(L^k)^{\boxtimes N_k}\longrightarrow M^{N_k}.$ Its squared norm is a scalar function on $M^{N_k}$, independent of the chosen orthonormal basis.

\begin{definition}
\label{def:Bergman-DPP}
The \emph{Bergman ensemble} associated with $(M,L^k,h^k,\nu_\psi)$ is the DPP corresponding to the probability measure on $M^{N_k}$ defined by
\[
    d\mathbb P_{k,\psi}(x_1,\ldots,x_{N_k})
    =
    \frac{1}{N_k!}
    \big\|
        \mathfrak S_{k,\psi}(x_1,\ldots,x_{N_k})
    \big\|_{h^k}^{2}
    \,d\nu_\psi(x_1)\cdots d\nu_\psi(x_{N_k}).
\]
\end{definition}

The density of the Bergman ensemble can also be rewritten with respect to the original measure $\mu$:
\[
    d\mathbb P_{k,\psi}(x_1,\ldots,x_{N_k})
    =
    \frac{1}{N_k!}
    \big\|
        \mathfrak S_{k,\psi}(x_1,\ldots,x_{N_k})
    \big\|_{h^k}^{2}
    e^{-\sum_{j=1}^{N_k}\psi(x_j)}
    \,d\mu(x_1)\cdots d\mu(x_{N_k}).
\]
By  Theorem~\ref{thm:intrinsic-determinantal-property}, the Bergman ensemble is a projection DPP, whose kernel is given by $B_{k,\psi}$. This process has a very natural physical motivation/interpretation \cite{DK,Kle,Kle2}: assume, for simplicity, that $M=\Sigma$ is a compact Riemann surface and that $L\to \Sigma$ is a positive Hermitian holomorphic line bundle. The curvature of $L^k$ represents a magnetic field of strength $k$. In the holomorphic quantization of a charged particle in this magnetic field, the lowest Landau level is identified with the finite-dimensional space $H^0(\Sigma,L^k),$ or, more generally, with $H^0(\Sigma,L^k\otimes K^s)$ if one includes a spin
coupling. Thus the Bergman projection is the orthogonal projection onto the lowest Landau level, and the Bergman kernel is the one-particle density matrix of the filled lowest Landau level.

The $N_k$-particle integer quantum Hall state at filling one is obtained by filling all one-particle states in the lowest Landau level. If $s_1,\ldots,s_{N_k}$ is an orthonormal basis of $H^0(\Sigma,L^k)$, the fermionic many-body wavefunction is the Slater determinant section~\eqref{eq:Slater}. Its squared norm, divided by $N_k!$, is exactly the density of the Bergman ensemble. Hence the Bergman ensemble is the position process of $N_k$ non-interacting fermions filling the lowest Landau level. Note that this determinantal description fails for general fractional quantum Hall effect.

This interpretation also explains why the kernel is naturally line-bundle-valued: local trivializations correspond to choices of gauge, and the scalar kernels written in a gauge transform by conjugation. The intrinsic determinant introduced above is the coordinate-free, gauge-invariant form of the fermionic correlation functions.

It is important to distinguish two equivalent descriptions of a change of weight. Suppose that $u:M\to\mathbb R$ is another bounded measurable function. Starting from the density of $\mathbb P_{k,\psi}$, one may form the tilted measure
\[
    d\mathbb Q(x_1,\ldots,x_{N_k})
    =
    \frac{1}{Z_{k,\psi}(u)}
    e^{-\sum_{j=1}^{N_k}u(x_j)}
    d\mathbb P_{k,\psi}(x_1,\ldots,x_{N_k}).
\]
Equivalently,
\[
    d\mathbb Q(x_1,\ldots,x_{N_k})
    =
    \frac{1}{\widetilde Z_{k,\psi}(u)}
    \big\|
        \mathfrak S_{k,\psi}(x_1,\ldots,x_{N_k})
    \big\|_{h^k}^{2}
    e^{-\sum_{j=1}^{N_k}(\psi+u)(x_j)}
    d\mu(x_1)\cdots d\mu(x_{N_k}).
\]

Although the basis used in the expression of $\mathfrak S_{k,\psi}$ is orthonormal for the old inner product $\langle\cdot,\cdot\rangle_{k,\psi}$, the squared Slater determinant defines the same alternating tensor up to a multiplicative constant when expressed in any basis of $H_k$. After normalization, the measure $\mathbb Q$ is therefore exactly the Bergman ensemble associated with the new inner product
\[
    \langle s,t\rangle_{k,\psi+u}
    =
    \int_M h_x^k(s(x),t(x))e^{-(\psi+u)(x)}\,d\mu(x).
\]

\begin{remark}\label{rmk:tilting-preserves-Bergman-DPP}
The tilted measure $\frac{1}{Z_{k,\psi}(u)}e^{-\sum_{j=1}^{N_k}u(x_j)}d\mathbb P_{k,\psi}$ is the Bergman ensemble $\mathbb P_{k,\psi+u}$. In particular, one-body changes of weight do not leave the class of projection DPPs; they only change the Hilbert structure and hence the Bergman projection kernel. This fact has a natural interpretation in the integer quantum Hall effect. In that setting, a change of magnetic potential, Hermitian metric, or background geometry changes the one-particle $L^2$-metric. The corresponding normalizing constant is the integer quantum Hall generating functional. In this language, a one-body perturbation of the many-body probability density is exactly the same as changing the Hermitian structure defining the Bergman projection. Thus the probabilistic tilting identity is the finite-$k$ version of the determinant-line variation which underlies the generating functional in the integer quantum Hall effect. In the large-$k$ regime, the asymptotics of this determinant-line norm are related to Quillen metrics, analytic torsion and anomaly formulas; see \cite{KMMW} for the relation between the integer quantum Hall generating functional, the Quillen metric and adiabatic curvature.
\end{remark}


\begin{proposition}\label{prop:palm-bergman}
Let $p_1,\ldots,p_\ell\in M$ be distinct points such that the evaluation map
\[
\operatorname{ev}_{\mathbf p}:H^0(M,L^k)\to L^k_{p_1}\oplus\cdots\oplus L^k_{p_\ell}
\]
has full rank. Let $H_k(-\mathbf p)=\{s\in H^0(M,L^k):s(p_1)=\cdots=s(p_\ell)=0\}.$ Then the reduced Palm process of the Bergman ensemble at $p_1,\ldots,p_\ell$ is the
projection DPP associated with the Hilbert space $H_k(-\mathbf p)$. Its kernel is the Schur complement
\[
B^{!,\mathbf p}_{k,\psi}(x,y)
=
B_{k,\psi}(x,y)
-
\sum_{a,b=1}^{\ell}
B_{k,\psi}(x,p_a)
\bigl(G_{k,\psi}(\mathbf p)^{-1}\bigr)_{ab}
B_{k,\psi}(p_b,y),
\]
where $G_{k,\psi}(\mathbf p)=\bigl(B_{k,\psi}(p_a,p_b)\bigr)_{1\leq a,b\leq \ell}.$
\end{proposition}

\begin{proof}
Write $H_k=H^0(M,L^k)$ with its $L^2(h^k,\nu_\psi)$-inner product, set $L^k_{\mathbf p}:=L^k_{p_1}\oplus\cdots\oplus L^k_{p_\ell},$ and let $\operatorname{ev}_{\mathbf p}:H_k\longrightarrow L^k_{\mathbf p}$ be the evaluation map. By assumption, $\operatorname{ev}_{\mathbf p}$ has full rank. The Gram operator at $\mathbf p$ is $G_{k,\psi}(\mathbf p)=\operatorname{ev}_{\mathbf p}\operatorname{ev}_{\mathbf p}^{*}\in \operatorname{End}(L^k_{\mathbf p}).$ Since $\operatorname{ev}_{\mathbf p}$ is surjective, $G_{k,\psi}(\mathbf p)$ is positive and invertible. In block form its $(a,b)$-entry is $B_{k,\psi}(p_a,p_b):L^k_{p_b}\longrightarrow L^k_{p_a}.$ Let $H_k(-\mathbf p)=\ker(\operatorname{ev}_{\mathbf p}).$ We first identify its orthogonal projection kernel. We have the orthogonal decomposition
\[
H_k=H_k(-\mathbf p)\oplus \operatorname{Ran}(\operatorname{ev}_{\mathbf p}^{*}).
\]
The orthogonal projection from $H_k$ onto
$\operatorname{Ran}(\operatorname{ev}_{\mathbf p}^{*})$ is
\[
P_{\mathbf p}
=
\operatorname{ev}_{\mathbf p}^{*}
G_{k,\psi}(\mathbf p)^{-1}
\operatorname{ev}_{\mathbf p}.
\]
Indeed, $P_{\mathbf p}$ is self-adjoint and idempotent, and its range is $\operatorname{Ran}(\operatorname{ev}_{\mathbf p}^{*})$. Hence the orthogonal projection onto $H_k(-\mathbf p)$ is
\[
\Pi_{k,\psi}^{-\mathbf p}
=
\Pi_{k,\psi}-P_{\mathbf p}.
\]
Its Schwartz kernel is therefore
\[
B_{k,\psi}^{!,\mathbf p}(x,y)
=
B_{k,\psi}(x,y)
-
\operatorname{ev}_{x}\operatorname{ev}_{\mathbf p}^{*}
G_{k,\psi}(\mathbf p)^{-1}
\operatorname{ev}_{\mathbf p}\operatorname{ev}_{y}^{*}.
\]
Writing this in block notation gives
\[
B_{k,\psi}^{!,\mathbf p}(x,y)
=
B_{k,\psi}(x,y)
-
\sum_{a,b=1}^{\ell}
B_{k,\psi}(x,p_a)
\bigl(G_{k,\psi}(\mathbf p)^{-1}\bigr)_{ab}
B_{k,\psi}(p_b,y).
\]
Thus the Schur-complement kernel is precisely the projection kernel of the Hilbert space $H_k(-\mathbf p)$.

It remains to identify this projection DPP with the reduced Palm process. By the general Palm formula for projection DPPs, recalled in Proposition~\ref{prop:palm-projection-DPP} and in the Schur-complement transfer principle of Section~\ref{sec:Palm-transfer}, the reduced Palm process at $\mathbf p$, provided
\[
\det_{\mathbf p}\bigl(B_{k,\psi}(p_a,p_b)\bigr)_{a,b=1}^{\ell}\neq0,
\]
is again a projection DPP, with kernel equal to the Schur complement of $G_{k,\psi}(\mathbf p)$ in the original kernel. The full-rank assumption on $\operatorname{ev}_{\mathbf p}$ is exactly the invertibility of this Gram operator.

Consequently, the reduced Palm process of the Bergman ensemble at $p_1,\ldots,p_\ell$, with those points removed, is the projection DPP associated with $H_k(-\mathbf p)$, and its kernel is the announced Schur-complement kernel.
\end{proof}

In the integer quantum Hall interpretation, Proposition~\ref{prop:palm-bergman} also has a natural physical meaning as a quasihole construction. Indeed, the filled lowest-Landau-level state is a Slater determinant. In the flat model at filling $\nu=1$, this determinant is the Vandermonde factor
\[
    \Psi_N(z_1,\ldots,z_N)
    =
    \prod_{i<j}(z_i-z_j)
\]
up to the Gaussian weight. If one conditions one particle to be at $\eta$, then
\[
    \Psi_N(z_1,\ldots,z_{N-1},\eta)
    =
    \prod_{i=1}^{N-1}(z_i-\eta)
    \prod_{1\le i<j\le N-1}(z_i-z_j),
\]
which is precisely the filled lowest-Landau-level state of the remaining particles with a quasihole inserted at $\eta$. The intrinsic version of this statement is that the reduced Palm process at $p$ is the projection DPP associated with $H_k(-p)$: the remaining one-particle states are constrained to vanish at $p$. Equivalently, the coherent state localized at $p$ has been removed from the filled lowest Landau level. The reduced Palm kernel
\[
    B_{k,\psi}^{!,p}(x,y)
    =
    B_{k,\psi}(x,y)
    -
    \frac{B_{k,\psi}(x,p)B_{k,\psi}(p,y)}{B_{k,\psi}(p,p)}
\]
is therefore the one-particle density matrix of this quasihole/Palm state.

\subsection{Limit theorems from Bergman kernel asymptotics}

We now explain how the abstract transfer principles of Section~\ref{sec:black-box-transfer} specialize to Bergman ensembles. On the analytic side, we quote standard inputs from Bergman-kernel asymptotics, Berezin--Toeplitz calculus and pluripotential theory. On the probabilistic side, the preceding sections show that each input has an immediate DPP consequence. The point is that the same finite-dimensional projection formalism simultaneously accounts for empirical measures, local correlations, Palm kernels, variances, cumulants and large deviations.

Throughout this subsection we impose the following smooth positive setting for simplicity. Let $M$ be a compact complex manifold of complex dimension $n$, let $(L,h)\to M$ be a positive holomorphic Hermitian line bundle, and let $\omega=c_1(L,h)$ be its curvature form, with the normalization fixed once and for all. Let $dV$ be a smooth positive volume form on $M$, let $\psi\in C^\infty(M,\mathbb R)$, and set $d\nu_\psi=e^{-\psi}dV.$ For every $k\geq 1$, let $B_{k,\psi}$ be the Bergman kernel of $H_k=H^0(M,L^k)\subset L^2(M,L^k;h^k,d\nu_\psi),$ and let $X_k$ be the associated Bergman ensemble. We write
\[
    N_k=\dim H^0(M,L^k),\qquad
    \widehat\mu_k=\frac1{N_k}\sum_{x\in X_k}\delta_x.
\]
By asymptotic Riemann--Roch \cite[Thm.~1.4.6]{MaMa},
\begin{equation}\label{eq:asympt-RR}
N_k=k^n\int_M \frac{\omega^n}{n!}+O(k^{n-1}).
\end{equation}
All constants below depend on the normalization convention for $c_1(L,h)$; the same convention is used in the Bergman kernel expansion and in the Riemann--Roch formula.

\subsubsection{Macroscopic limit}


\begin{theorem}[{\cite[Thms.~4.1.1--4.1.2]{MaMa}}]
\label{thm:MM-diagonal}
Let $X$ be a compact K\"ahler manifold, let $(L,h^L)\to X$ be a positive Hermitian holomorphic line bundle, and set
\[
    \omega=c_1(L,h^L)=\frac{\sqrt{-1}}{2\pi}R^L.
\]
Let $(E,h^E)\to X$ be a Hermitian holomorphic vector bundle. Let $P_k^E(x,y)$ be the Bergman kernel of the orthogonal projection onto $H^0(X,L^k\otimes E)$, computed with respect to $dV_\omega=\omega^n/n!$. Then, for every $m\geq 0$,
\[
    P_k^E(x,x)
    =
    k^n b_0(x)+k^{n-1}b_1(x)+\cdots+k^{n-m}b_m(x)
    +O(k^{n-m-1})
\]
in $C^\infty(X,\operatorname{End}E)$, with $b_0(x)=\operatorname{Id}_{E_x}.$
\end{theorem}

\begin{theorem}\label{thm:bergman-macroscopic-convergence}
In the smooth positive setting, the empirical measures of the Bergman ensemble
satisfy
\[
    \widehat\mu_k
    =
    \frac1{N_k}\sum_{x\in X_k}\delta_x
    \Longrightarrow
    \frac{\omega^n}{\int_M\omega^n}
\]
in probability.
\end{theorem}

\begin{proof}
Apply Theorem~\ref{thm:MM-diagonal} with the trivial holomorphic line bundle $E=M\times\mathbb C$ endowed with the metric
\[
    |1|_{h^E}^2=\frac{d\nu_\psi}{dV_\omega}.
\]
If $P_k^E$ denotes the Bergman kernel of $L^k\otimes E$ with respect to $dV_\omega$, then
\[
    P_k^E(x,x)\,dV_\omega(x)
    =
    B_{k,\psi}(x,x)\,d\nu_\psi(x).
\]
The diagonal expansion therefore gives
\[
    B_{k,\psi}(x,x)d\nu_\psi(x)
    =
    k^n\frac{\omega^n}{n!}+O(k^{n-1})
\]
as smooth densities. Together with~\eqref{eq:asympt-RR}, we get
\[
    \frac1{N_k}B_{k,\psi}(x,x)d\nu_\psi(x)
    \Longrightarrow
    \frac{\omega^n}{\int_M\omega^n}.
\]
The result follows from Theorem~\ref{thm:diagonal-transfer}.
\end{proof}



\subsubsection{Variance bounds}

For a bundle-valued kernel we write
\[
    |B_{k,\psi}(x,y)|^2
    =
    \operatorname{Tr}\bigl(B_{k,\psi}(x,y)B_{k,\psi}(y,x)\bigr).
\]
Since $L_x^k$ is one-dimensional, this is simply the squared Hilbert--Schmidt norm of the map $L_y^k\to L_x^k$.

\begin{proposition}\label{prop:quadratic-localization}
In the smooth positive setting, for every Riemannian distance $d$ on $M$ and every integer $r\geq 0$,
\[
    \int_{M\times M}
    d(x,y)^r |B_{k,\psi}(x,y)|^2
    d\nu_\psi(x)d\nu_\psi(y)
    =
    O(k^{-r/2}N_k).
\]
\end{proposition}

\begin{proof}
As in the proof of Theorem~\ref{thm:bergman-macroscopic-convergence}, the smooth density $d\nu_\psi$ may be absorbed into a trivial auxiliary bundle, so we may use the standard off-diagonal estimates for Bergman kernels. The near-diagonal expansion and the off-diagonal estimates of Ma--Marinescu 
\cite[Thms.~4.1.24 and~4.2.1]{MaMa} give constants $C,c>0$ such that, for $x,y$ sufficiently close and all $k$,
\[
    |B_{k,\psi}(x,y)|
    \leq
    Ck^n\exp(-ckd(x,y)^2),
\]
while outside a fixed small neighborhood of the diagonal the kernel is $O(k^{-N})$ for every $N$. Hence
\[
\begin{aligned}
&\int_{M\times M}
    d(x,y)^r |B_{k,\psi}(x,y)|^2
    d\nu_\psi(x)d\nu_\psi(y) \\
&\qquad \leq
    Ck^{2n}\int_M\int_M d(x,y)^r e^{-ckd(x,y)^2}\,dV(y)dV(x)
    +O(k^{-N})
\end{aligned}
\]
for every $N$. In normal coordinates around $x$, the inner integral is
\[
    O\left(k^{-(n+r/2)}\right),
\]
because $M$ has real dimension $2n$. Thus the whole expression is $O(k^{n-r/2})$. Since $N_k=k^n\int_M\omega^n/n!+O(k^{n-1})$, this is $O(k^{-r/2}N_k)$.
\end{proof}

\begin{proposition}\label{prop:bergman-variance-lipschitz}
Let $f:M\to\mathbb R$ be Lipschitz. Then
\[
    \operatorname{Var}(\widehat\mu_k(f))
    \leq
    \frac{C\operatorname{Lip}(f)^2}{kN_k}
    =
    O(k^{-n-1}).
\]
\end{proposition}

\begin{proof}
By the projection-DPP variance identity,
\[
    \operatorname{Var}\left(\widehat\mu_k(f)\right)
    =
    \frac1{2N_k^2}
    \int_{M\times M}
        |f(x)-f(y)|^2 |B_{k,\psi}(x,y)|^2
        d\nu_\psi(x)d\nu_\psi(y).
\]
Since $f$ is Lipschitz,
\[
    |f(x)-f(y)|^2\leq \operatorname{Lip}(f)^2 d(x,y)^2.
\]
The result follows from Proposition~\ref{prop:quadratic-localization} with $r=2$.
\end{proof}

\subsubsection{Local universality and local fluctuation expansions}

The next result is the smooth positive version of the near-diagonal Bergman kernel expansion. We state it in the form needed for the determinantal process: the kernel is written with respect to Lebesgue measure in a local chart.

Let $p\in M$. Choose holomorphic coordinates $z=(z_1,\ldots,z_n)$ centered at $p$ and a local holomorphic frame $e$ of $L$ such that
\[
    |e(z)|_h^2=e^{-\varphi(z)},
    \qquad
    \varphi(z)=H_p(z,z)+O(|z|^3),
\]
where $H_p$ is a positive Hermitian form on $\mathbb C^n$. Write
\[
    H_p(u,v)=\sum_{\alpha,\beta}H_{\alpha\bar\beta}(p)
    u_\alpha\overline{v_\beta}.
\]
Write $dV=a(z)\,d\lambda(z)$ in the chosen coordinate chart, where $d\lambda$ denotes Lebesgue measure on $\mathbb C^n$. If $s_\ell^{(k)}=f_\ell^{(k)}e^{\otimes k}$ is an orthonormal basis of $H^0(M,L^k)$, set
\[
\mathcal K_{k,\psi}^{\mathrm{Leb}}(z,w)
=
a(z)^{1/2}e^{-\psi(z)/2}e^{-k\varphi(z)/2}
\left(
\sum_{\ell=1}^{N_k}
f_\ell^{(k)}(z)\overline{f_\ell^{(k)}(w)}
\right)
e^{-k\varphi(w)/2}e^{-\psi(w)/2}a(w)^{1/2}.
\]
This is the local scalar projection kernel written with respect to Lebesgue measure.

\begin{theorem}[{\cite[Thms.~4.1.24 and~4.2.1]{MaMa}}]\label{thm:MM-near-diagonal}
There exist smooth kernels
\[
    K_{p,j}(u,v),\qquad j\geq 0,
\]
defined locally on $\mathbb C^n\times\mathbb C^n$, such that, for every $N\geq 0$,
\[
    k^{-n}
    \mathcal K_{k,\psi}^{\mathrm{Leb}}
    \left(\frac{u}{\sqrt k},\frac{v}{\sqrt k}\right)
    =
    \sum_{j=0}^N k^{-j/2}K_{p,j}(u,v)
    +O(k^{-(N+1)/2})
\]
locally uniformly, with all derivatives, on compact subsets of
$\mathbb C^n\times\mathbb C^n$. The leading coefficient is
\[
    K_{p,0}(u,v)=K_p^{\mathrm{BF}}(u,v),
\]
where
\[
    K_p^{\mathrm{BF}}(u,v)
    =
    \frac{\det H_p}{\pi^n}
    \exp\left(
        H_p(u,v)
        -\frac12H_p(u,u)
        -\frac12H_p(v,v)
    \right).
\]
\end{theorem}

\begin{theorem}\label{thm:BF-local-universality}
For every fixed $m\geq 1$, there exist smooth functions
\[
    R_{m,j,p}:(\mathbb C^n)^m\longrightarrow \mathbb C,
    \qquad j\geq 0,
\]
such that, for every $N\geq 0$,
\[
    k^{-nm}\rho^{\mathrm{loc}}_{m,k}
    \left(\frac{u_1}{\sqrt k},\ldots,\frac{u_m}{\sqrt k}\right)
    =
    \sum_{j=0}^N k^{-j/2}
    R_{m,j,p}(u_1,\ldots,u_m)
    +
    O(k^{-(N+1)/2})
\]
locally uniformly on compact subsets of $(\mathbb C^n)^m$.

The coefficients are obtained by taking determinants of the formal kernel expansion:
\[
    R_{m,j,p}(u_1,\ldots,u_m)
    =
    [t^j]\,
    \det\left(
        K_{p,0}(u_a,u_b)
        +
        \sum_{q\geq 1}t^q K_{p,q}(u_a,u_b)
    \right)_{1\leq a,b\leq m}.
\]
\end{theorem}

\begin{proof}
By the determinantal formula,
\[
    \rho^{\mathrm{loc}}_{m,k}(z_1,\ldots,z_m)
    =
    \det\left(
        \mathcal K_{k,\psi}^{\mathrm{Leb}}(z_a,z_b)
    \right)_{a,b=1}^m .
\]
Taking $z_a=u_a/\sqrt k$ gives
\[
    k^{-nm}\rho^{\mathrm{loc}}_{m,k}
    \left(\frac{u_1}{\sqrt k},\ldots,\frac{u_m}{\sqrt k}\right)
    =
    \det\left(
        k^{-n}\mathcal K_{k,\psi}^{\mathrm{Leb}}
        \left(\frac{u_a}{\sqrt k},\frac{u_b}{\sqrt k}\right)
    \right)_{a,b=1}^m .
\]
Theorem~\ref{thm:MM-near-diagonal} gives a full expansion of each entry of this matrix. Since the determinant is a polynomial in the entries, the determinant also admits a full expansion, and its coefficients are exactly the stated formal determinant coefficients.
\end{proof}

Note that Theorem~\ref{thm:BF-local-universality} gives in particular the usual Bargmann--Fock universality
\[
    R_{m,0,p}(u_1,\ldots,u_m)
    =
    \det\left(K^{\mathrm{BF}}_p(u_a,u_b)\right)_{a,b=1}^m.
\]
It is useful to spell out one immediate probabilistic consequence of the expansion. Let
\[
    \Xi_{k,p}
    =
    \sum_{x\in X_k\cap U}\delta_{\sqrt k\,z(x)}
\]
be the rescaled local point process in the chart. If $\chi\in C_c^\infty(\mathbb C^n)$ and $k$ is large enough so that $\operatorname{supp}\chi$ is contained in the rescaled coordinate domain, define the local linear statistic
\[
    S_{k,p}(\chi)
    =
    \int_{\mathbb C^n}\chi(u)\,d\Xi_{k,p}(u)
    =
    \sum_{x\in X_k}\chi(\sqrt k\,z(x)).
\]
Let $\chi_1,\ldots,\chi_r\in C_c^\infty(\mathbb C^n)$. For every $r\geq 1$, the joint factorial moments
\[
    \mathbb E\left[
        \sum_{u_1,\ldots,u_r\in \Xi_{k,p}}^{\neq}
        \chi_1(u_1)\cdots \chi_r(u_r)
    \right]
\]
admit complete asymptotic expansions in descending half-integer powers of $k$. More precisely, for every $N\geq 0$,
\begin{align*}
\mathbb E & \left[\sum_{u_1,\ldots,u_r\in \Xi_{k,p}}^{\neq} \chi_1(u_1)\cdots \chi_r(u_r)\right]\\
& =\sum_{j=0}^N k^{-j/2}\int_{(\mathbb C^n)^r}\chi_1(u_1)\cdots\chi_r(u_r)R_{r,j,p}(u_1,\ldots,u_r)\,d\lambda(u_1)\cdots d\lambda(u_r)+O(k^{-(N+1)/2}).
\end{align*}
Consequently, all joint moments and all joint cumulants of $S_{k,p}(\chi_1),\ldots,S_{k,p}(\chi_r)$ admit complete asymptotic expansions in descending half-integer powers of $k$. Their leading terms are the corresponding moments and cumulants of the Bargmann--Fock DPP with kernel $K_p^{\mathrm{BF}}$.

In particular, for real-valued $\chi_1,\chi_2\in C_c^\infty(\mathbb C^n)$,
\[
    \operatorname{Cov}\bigl(S_{k,p}(\chi_1),S_{k,p}(\chi_2)\bigr)
    =
    \operatorname{Cov}_{\mathrm{BF},p}\bigl(S(\chi_1),S(\chi_2)\bigr)
    +O(k^{-1/2}),
\]
and in fact the left-hand side has a full expansion in powers of $k^{-1/2}$. The leading covariance is
\[
\begin{aligned}
    \operatorname{Cov}_{\mathrm{BF},p}\bigl(S(\chi_1),S(\chi_2)\bigr)
    &=
    \int_{\mathbb C^n}
        \chi_1(u)\chi_2(u)K_p^{\mathrm{BF}}(u,u)\,d\lambda(u)  \\
    &\quad
    -
    \int_{\mathbb C^n\times\mathbb C^n}
        \chi_1(u)\chi_2(v)
        \left|K_p^{\mathrm{BF}}(u,v)\right|^2
        \,d\lambda(u)d\lambda(v).
\end{aligned}
\]
Equivalently, if $\chi_1=\chi_2=\chi$, then
\begin{equation}
\operatorname{Var}_{\mathrm{BF},p}(S(\chi))=\frac12\int_{\mathbb C^n\times\mathbb C^n}|\chi(u)-\chi(v)|^2\left|K_p^{\mathrm{BF}}(u,v)\right|^2\,d\lambda(u)d\lambda(v).
\end{equation}
Note that Theorem~\ref{thm:BF-local-universality}, combined with the previous near-diagonal expansion and off-diagonal estimates used in Proposition~\ref{prop:quadratic-localization}, implies the required local trace-norm convergence on compact subsets; hence Proposition~\ref{prop:local-process-transfer} applies and we obtain a convergence of the rescaled DPPs. We have even better than that: the convergence of all reduced Palm processes.

\begin{theorem}\label{thm:Palm-BF-universality}
Let $q_1,\ldots,q_\ell\in\mathbb C^n$ be distinct points and set $p_{a,k}=z^{-1}(q_a/\sqrt k).$ Let $X_k^{!,p_k}$ be the reduced Palm process of the Bergman ensemble conditioned to contain $p_{1,k},\ldots,p_{\ell,k}$, with these points removed. Then, the local correlation functions of the reduced Palm Bergman process converge to those of the reduced Palm Bargmann--Fock process.
\end{theorem}

\begin{proof}
It boils down to prove
\[
    k^{-n}\mathcal K_{k,\psi}^{\mathrm{Leb},!,p_k}
    \left(
        \frac{u}{\sqrt k},\frac{v}{\sqrt k}
    \right)
    \longrightarrow
    K^{\mathrm{BF},!,q}_p(u,v)
\]
locally uniformly, where
\[
    K^{\mathrm{BF},!,q}_p(u,v)
    =
    K^{\mathrm{BF}}_p(u,v)
    -
    \sum_{a,b=1}^\ell
    K^{\mathrm{BF}}_p(u,q_a)
    \left(G_p(q)^{-1}\right)_{ab}
    K^{\mathrm{BF}}_p(q_b,v),
\]
and
\[
    G_p(q)=
    \left(K^{\mathrm{BF}}_p(q_a,q_b)\right)_{1\leq a,b\leq \ell}.
\]
The Bargmann--Fock kernel is strictly positive definite, hence $G_p(q)$ is invertible for distinct $q_1,\ldots,q_\ell$. The reduced Palm kernel is the Schur complement of the finite Gram matrix at the conditioned points, by Proposition~\ref{prop:palm-bergman}. The result follows by applying Proposition~\ref{prop:palm-transfer} to the leading kernel convergence contained in Theorem~\ref{thm:MM-near-diagonal}.
\end{proof}

\subsubsection{Cumulants}

We now turn from correlation functions to linear statistics. The relevant analytic object is no longer only the Bergman kernel but the Berezin--Toeplitz calculus. The finite-dimensional DPP identity expresses cumulants as traces of products of compressed multiplication operators. Toeplitz asymptotics then convert these exact trace formulas into semiclassical expansions.

\begin{theorem}[\cite{MaMa2}]
\label{thm:MM-toeplitz-product}
Let $f,g\in C^\infty(M)$. Then
\[
    T_{f,k,\psi}T_{g,k,\psi}
    \sim
    \sum_{j\geq 0}k^{-j}T_{C_j(f,g),k,\psi},
\]
where $C_j$ are bidifferential operators and $C_0(f,g)=fg.$ More generally, products of finitely many Toeplitz operators admit a complete asymptotic expansion, whose leading symbol is the product of the symbols.
\end{theorem}

\begin{corollary}\label{cor:toeplitz-trace-asymp}
For $f_1,\ldots,f_r\in C^\infty(M)$,
\[
    \operatorname{Tr}(T_{f_1,k,\psi}\cdots T_{f_r,k,\psi})
    \sim
    k^n\sum_{j\geq 0}k^{-j}
    \int_M b_j(f_1,\ldots,f_r)\frac{\omega^n}{n!},
\]
with
\[
    b_0(f_1,\ldots,f_r)=f_1\cdots f_r.
\]
\end{corollary}

\begin{proposition}\label{prop:cumulant-asymptotics}
Let $f\in C^\infty(M,\mathbb R)$, and set $S_{f,k}=\sum_{x\in X_k}f(x).$ For every $\ell\geq 1$, the cumulant $\kappa_{\ell,k}(f)$ admits a full asymptotic expansion in descending powers of $k$. Moreover,
\[
    \kappa_{1,k}(f)
    =
    k^n\int_M f\,\frac{\omega^n}{n!}
    +O(k^{n-1}),
\]
whereas, for every $\ell\geq 2$, $\kappa_{\ell,k}(f)=O(k^{n-1}).$
\end{proposition}

\begin{proof}
By Proposition~\ref{prop:cumulant-trace-formula} applied to the Bergman projection,
\[
    \kappa_{\ell,k}(f)
    =
    \sum_{r=1}^{\ell}
    \frac{(-1)^{r-1}}{r}
    \sum_{\substack{m_1+\cdots+m_r=\ell\\ m_i\geq 1}}
    \frac{\ell!}{m_1!\cdots m_r!}
    \operatorname{Tr}
    \bigl(
        T_{f^{m_1},k,\psi}\cdots T_{f^{m_r},k,\psi}
    \bigr).
\]
For each fixed composition $m_1+\cdots+m_r=\ell$, Corollary~\ref{cor:toeplitz-trace-asymp} gives a full asymptotic expansion
\[
    \operatorname{Tr}
    \bigl(
        T_{f^{m_1},k,\psi}\cdots T_{f^{m_r},k,\psi}
    \bigr)
    \sim
    k^n\sum_{j\geq 0}k^{-j}
    \int_M
    b_j(f^{m_1},\ldots,f^{m_r})\frac{\omega^n}{n!}.
\]
Since only finitely many compositions occur for fixed $\ell$, this gives a full asymptotic expansion for $\kappa_{\ell,k}(f)$.

It remains to identify the coefficient of $k^n$. Since $b_0(g_1,\ldots,g_r)=g_1\cdots g_r,$ the leading contribution equals
\[
    k^n\int_M f^\ell\,\frac{\omega^n}{n!}\,A_\ell,
\]
where
\[
    A_\ell
    =
    \sum_{r=1}^{\ell}
    \frac{(-1)^{r-1}}{r}
    \sum_{\substack{m_1+\cdots+m_r=\ell\\ m_i\geq 1}}
    \frac{\ell!}{m_1!\cdots m_r!}.
\]
The number $A_\ell/\ell!$ is the coefficient of $t^\ell$ in
\[
    \sum_{r\geq 1}
    \frac{(-1)^{r-1}}{r}(e^t-1)^r
    =
    \log(1+e^t-1)
    =
    t.
\]
Hence $A_1=1$ and $A_\ell=0$ for every $\ell\geq 2$. This proves the claimed leading orders.
\end{proof}

\begin{remark}
The preceding proof gives more than the order estimates. If
\[
    \operatorname{Tr}
    \bigl(T_{g_1,k,\psi}\cdots T_{g_r,k,\psi}\bigr)
    \sim
    k^n\sum_{j\geq 0}k^{-j}
    \int_M b_j(g_1,\ldots,g_r)\frac{\omega^n}{n!},
\]
then the coefficient of $k^{n-j}$ in $\kappa_{\ell,k}(f)$ is
\[
    \sum_{r=1}^{\ell}
    \frac{(-1)^{r-1}}{r}
    \sum_{\substack{m_1+\cdots+m_r=\ell\\ m_i\geq 1}}
    \frac{\ell!}{m_1!\cdots m_r!}
    \int_M
    b_j(f^{m_1},\ldots,f^{m_r})\frac{\omega^n}{n!}.
\]
For instance,
\[
    \kappa_{2,k}(f)
    =
    \operatorname{Tr}(T_{f^2,k,\psi})
    -
    \operatorname{Tr}(T_{f,k,\psi}^2),
\]
so its leading $k^n$-term cancels. This cancellation is the Toeplitz form of the rigidity of projection DPPs: the macroscopic statistic is deterministic at principal semiclassical order, and fluctuations only appear at subprincipal order.
\end{remark}

This cumulant expansion should be viewed as a smooth positive counterpart of the fluctuation results proved by Berman in greater generality. The present argument does not aim to optimize regularity assumptions; it shows that, once Toeplitz trace expansions are available, the probabilistic cumulants follow formally.

\subsubsection{Large deviations}

The large-deviation regime is different from the preceding smooth local asymptotics. The relevant exponential tilts are of the form
\[
    \exp\left(k\sum_{i=1}^{N_k} f(x_i)\right),
\]
and therefore change the geometric weight $k\varphi$ at leading order. The correct analytic input is no longer the local Bergman kernel expansion, but the asymptotic behavior of determinant functionals, or equivalently of volumes of balls of holomorphic sections.

We state the input in a form adapted to the transfer theorem. Fix a reference measure $\mu$ satisfying the Bernstein--Markov property in the sense of \cite{BB}, and write $Z_k(\varphi)$ for the unnormalized partition function
\[
    Z_k(\varphi)
    =
    \int_{M^{N_k}}
    \|\mathfrak S(x_1,\ldots,x_{N_k})\|^2
    \exp\left(-k\sum_{i=1}^{N_k}\varphi(x_i)\right)
    d\mu(x_1)\cdots d\mu(x_{N_k}),
\]
where $\mathfrak S$ is formed from any fixed basis of $H^0(M,L^k)$. The ratio $Z_k(\varphi-f)/Z_k(\varphi)$ is independent of this fixed basis.

\begin{theorem}[{\cite{BB}}]
\label{thm:BB-determinant}
Let $L$ be big over a compact complex manifold $M$, let $\varphi\in C(M)$ be a continuous weight, and let $\mu$ be a probability measure on $M$ satisfying the Bernstein--Markov property with respect to $(L,\varphi)$. Then, for every continuous real function $f$,
\[
\lim_{k\to\infty}\frac{1}{kN_k}\log \frac{Z_k(\varphi-f)}{Z_k(\varphi)}=\mathcal E_{\mathrm{eq}}(\varphi)-\mathcal E_{\mathrm{eq}}(\varphi-f),
\]
where $\mathcal E_{\mathrm{eq}}$ is normalized so that $\left.\frac{d}{dt}\right|_{t=0}\mathcal E_{\mathrm{eq}}(\varphi+tv)=\int_M v\,d\mu_{\mathrm{eq}}(\varphi).$
\end{theorem}

\begin{corollary}\label{cor:Bergman-LDP-transfer}
Under the hypotheses of Theorem~\ref{thm:BB-determinant}, the empirical measures of the Bergman ensemble with weight $\varphi$ satisfy a large-deviation principle at speed $kN_k$ with good rate function
\[
    I_\varphi(\nu)=\sup_{f\in C(M,\mathbb R)}\left\{\int_M f\,d\nu-\bigl(\mathcal E_{\mathrm{eq}}(\varphi)-\mathcal E_{\mathrm{eq}}(\varphi-f)\bigr)\right\}.
\]
The rate function is minimized at the equilibrium measure $\mu_{\mathrm{eq}}(\varphi)$. In the smooth strictly positive case, this equilibrium measure is the normalized curvature volume.
\end{corollary}

\begin{proof}
For the Bergman ensemble with weight $\varphi$,
\[
    \mathbb E_\varphi
    \left[
        \exp\left(k\sum_{i=1}^{N_k}f(x_i)\right)
    \right]
    =
    \frac{Z_k(\varphi-f)}{Z_k(\varphi)}.
\]
Since $k\sum_i f(x_i)=kN_k\,\widehat\mu_k(f)$, Theorem~\ref{thm:BB-determinant} identifies the limiting log-Laplace functional in Theorem~\ref{thm:ldp-transfer} with $\Lambda_\varphi(f)=\mathcal E_{\mathrm{eq}}(\varphi)-\mathcal E_{\mathrm{eq}}(\varphi-f).$ The large-deviation principle follows from Theorem~\ref{thm:ldp-transfer}. The differentiability formula for $\mathcal E_{\mathrm{eq}}$ gives $D\Lambda_\varphi(0)(f)=\int_M f\,d\mu_{\mathrm{eq}}(\varphi),$ which identifies the minimizer of the Legendre transform. In the smooth positive setting the equilibrium envelope is the weight itself, so $\mu_{\mathrm{eq}}(\varphi)$ is the normalized Monge--Amp\`ere, equivalently curvature-volume, measure.
\end{proof}

\begin{remark}
Berman proved the corresponding large-deviation principle in a substantially more general pluripotential framework \cite{Ber6}. The point of the present formulation is not to replace that theorem, but to show exactly how the determinant asymptotic enters the finite-dimensional DPP formalism: after the identity
\[
    \mathbb E_\varphi\left[e^{kN_k\widehat\mu_k(f)}\right]
    =
    Z_k(\varphi-f)/Z_k(\varphi),
\]
the probabilistic part is precisely the abstract transfer theorem.
\end{remark}

\begin{remark}
A bounded tilt $e^{-\sum_i u(x_i)}$ changes the auxiliary weight $\psi$ to $\psi+u$, as in Remark~\ref{rmk:tilting-preserves-Bergman-DPP}. By contrast, the large-deviation tilt changes the geometric weight $k\varphi$ to $k(\varphi+V)$. Thus the latter belongs to the pluripotential regime. If $\varphi+V$ is strictly positive, the limiting equilibrium measure is the curvature volume of $\varphi+V$. For merely continuous $V$, it is the Monge--Amp\`ere measure of the corresponding equilibrium envelope.
\end{remark}

\section*{Funding}

This work was supported by the Agence Nationale de la Recherche
[grant number ANR-20-CHIA-0002].

\section*{Acknowledgments}

I thank Semyon Klevtsov for having introduced me to the quantum Hall effect, and for many stimulating discussions, as the first version of this paper was written during his project ``Geometry of quantum Hall states'' supported by the IdEx program and the USIAS fellowship of Universit\'e de Strasbourg. I also thank Michele Ancona and Yohann Le Floch for discussions about probability and complex geometry, and Rapha\"el Butez for several discussions about random matrices and random polynomials. The later version of this paper also benefited from stimulating discussions with J\'er\'emie Bouttier, Laurent Charles, Benoit Estienne, David Garc\'ia-Zelada, Louis Ioos, Thierry L\'evy and Elias Nohra.

\bibliographystyle{alpha}
\bibliography{DPP_Bergman}

@online{DPP-Fermion,
      title={From point processes to quantum optics and back}, 
      author={Rémi Bardenet and Alexandre Feller and Jérémie Bouttier and Pascal Degiovanni and Adrien Hardy and Adam Rançon and Benjamin Roussel and Grégory Schehr and Christoph I. Westbrook},
      year={2022},
      eprint={2210.05522},
      note={arxiv:2210.05522},
      archivePrefix={arXiv},
      primaryClass={math-ph},
      url={https://arxiv.org/abs/2210.05522}, 
}

@article {BB,
    AUTHOR = {Berman, Robert and Boucksom, S\'{e}bastien},
     TITLE = {Growth of balls of holomorphic sections and energy at
              equilibrium},
   JOURNAL = {Invent. Math.},
  FJOURNAL = {Inventiones Mathematicae},
    VOLUME = {181},
      YEAR = {2010},
    NUMBER = {2},
     PAGES = {337--394},
      ISSN = {0020-9910},
   MRCLASS = {32L05 (32L10 32U15 32W20 58J52)},
  MRNUMBER = {2657428},
MRREVIEWER = {Norman Levenberg},
       DOI = {10.1007/s00222-010-0248-9},
       URL = {https://doi.org/10.1007/s00222-010-0248-9},
}

@article {Ber6,
    AUTHOR = {Berman, Robert J.},
     TITLE = {Determinantal point processes and fermions on complex
              manifolds: large deviations and bosonization},
   JOURNAL = {Comm. Math. Phys.},
  FJOURNAL = {Communications in Mathematical Physics},
    VOLUME = {327},
      YEAR = {2014},
    NUMBER = {1},
     PAGES = {1--47},
      ISSN = {0010-3616},
   MRCLASS = {81V35},
  MRNUMBER = {3177931},
       DOI = {10.1007/s00220-014-1891-6},
       URL = {https://doi.org/10.1007/s00220-014-1891-6},
}

@incollection {Ber7,
    AUTHOR = {Berman, Robert J.},
     TITLE = {Determinantal point processes and fermions on polarized
              complex manifolds: bulk universality},
 BOOKTITLE = {Algebraic and analytic microlocal analysis},
    SERIES = {Springer Proc. Math. Stat.},
    VOLUME = {269},
     PAGES = {341--393},
 PUBLISHER = {Springer, Cham},
      YEAR = {2018},
   MRCLASS = {32A25 (32W20 60B20 60G55)},
  MRNUMBER = {3903320},
MRREVIEWER = {Daniel Belti\c{t}\u{a}},
       DOI = {10.1007/978-3-030-01588-6\_5},
       URL = {https://doi.org/10.1007/978-3-030-01588-6_5},
}

@article {BS,
    AUTHOR = {Borodin, A. and Soshnikov, A.},
     TITLE = {Janossy densities. {I}. {D}eterminantal ensembles},
   JOURNAL = {J. Statist. Phys.},
  FJOURNAL = {Journal of Statistical Physics},
    VOLUME = {113},
      YEAR = {2003},
    NUMBER = {3-4},
     PAGES = {595--610},
      ISSN = {0022-4715},
   MRCLASS = {60G55 (35Q15 82B44)},
  MRNUMBER = {2013698},
MRREVIEWER = {Wolfgang Freudenberg},
       DOI = {10.1023/A:1026025003309},
       URL = {https://doi.org/10.1023/A:1026025003309},
}

@article {Buf23,
    AUTHOR = {Bufetov, Alexander I.},
     TITLE = {The conditional measures for the determinantal point process
              with the {B}ergman kernel},
   JOURNAL = {Atti Accad. Naz. Lincei Rend. Lincei Mat. Appl.},
  FJOURNAL = {Atti della Accademia Nazionale dei Lincei. Rendiconti Lincei.
              Matematica e Applicazioni},
    VOLUME = {34},
      YEAR = {2023},
    NUMBER = {1},
     PAGES = {159--173},
      ISSN = {1120-6330,1720-0768},
   MRCLASS = {60G57 (30H20 60G55)},
  MRNUMBER = {4628738},
       DOI = {10.4171/rlm/1002},
       URL = {https://doi.org/10.4171/rlm/1002},
}

@article {BQ,
    AUTHOR = {Bufetov, Alexander I. and Qiu, Yanqi},
     TITLE = {Determinantal point processes associated with {H}ilbert spaces
              of holomorphic functions},
   JOURNAL = {Comm. Math. Phys.},
  FJOURNAL = {Communications in Mathematical Physics},
    VOLUME = {351},
      YEAR = {2017},
    NUMBER = {1},
     PAGES = {1--44},
      ISSN = {0010-3616},
   MRCLASS = {60G55 (30H20 60G57)},
  MRNUMBER = {3613499},
       DOI = {10.1007/s00220-017-2840-y},
       URL = {https://doi.org/10.1007/s00220-017-2840-y},
}

@article {BFQ,
    AUTHOR = {Bufetov, Alexander I. and Fan, Shilei and Qiu, Yanqi},
     TITLE = {Equivalence of {P}alm measures for determinantal point
              processes governed by {B}ergman kernels},
   JOURNAL = {Probab. Theory Related Fields},
  FJOURNAL = {Probability Theory and Related Fields},
    VOLUME = {172},
      YEAR = {2018},
    NUMBER = {1-2},
     PAGES = {31--69},
      ISSN = {0178-8051,1432-2064},
   MRCLASS = {60G55 (32A36)},
  MRNUMBER = {3851829},
       DOI = {10.1007/s00440-017-0803-z},
       URL = {https://doi.org/10.1007/s00440-017-0803-z},
}

@article {Bor,
    AUTHOR = {Borodin, Alexei},
     TITLE = {Biorthogonal ensembles},
   JOURNAL = {Nuclear Phys. B},
  FJOURNAL = {Nuclear Physics. B. Theoretical, Phenomenological, and
              Experimental High Energy Physics. Quantum Field Theory and
              Statistical Systems},
    VOLUME = {536},
      YEAR = {1999},
    NUMBER = {3},
     PAGES = {704--732},
      ISSN = {0550-3213},
   MRCLASS = {82B41 (15A52)},
  MRNUMBER = {1663328},
MRREVIEWER = {Oleksiy Khorunzhiy},
       DOI = {10.1016/S0550-3213(98)00642-7},
       URL = {https://doi.org/10.1016/S0550-3213(98)00642-7},
}

@book {DZ,
    AUTHOR = {Dembo, Amir and Zeitouni, Ofer},
     TITLE = {Large deviations techniques and applications},
    SERIES = {Applications of Mathematics (New York)},
    VOLUME = {38},
   EDITION = {Second},
 PUBLISHER = {Springer-Verlag, New York},
      YEAR = {1998},
     PAGES = {xvi+396},
      ISBN = {0-387-98406-2},
   MRCLASS = {60F10},
  MRNUMBER = {1619036},
       DOI = {10.1007/978-1-4612-5320-4},
       URL = {https://doi.org/10.1007/978-1-4612-5320-4},
}

@article {DK,
    AUTHOR = {Douglas, Michael R. and Klevtsov, Semyon},
     TITLE = {Bergman kernel from path integral},
   JOURNAL = {Comm. Math. Phys.},
  FJOURNAL = {Communications in Mathematical Physics},
    VOLUME = {293},
      YEAR = {2010},
    NUMBER = {1},
     PAGES = {205--230},
      ISSN = {0010-3616},
   MRCLASS = {32L10 (32Q15 47B35 58J37 81Q60 81T20)},
  MRNUMBER = {2563804},
MRREVIEWER = {Valentino Tosatti},
       DOI = {10.1007/s00220-009-0915-0},
       URL = {https://doi.org/10.1007/s00220-009-0915-0},
}

@book {Dei,
    AUTHOR = {Deift, P. A.},
     TITLE = {Orthogonal polynomials and random matrices: a
              {R}iemann-{H}ilbert approach},
    SERIES = {Courant Lecture Notes in Mathematics},
    VOLUME = {3},
 PUBLISHER = {New York University, Courant Institute of Mathematical
              Sciences, New York; American Mathematical Society, Providence,
              RI},
      YEAR = {1999},
     PAGES = {viii+273},
      ISBN = {0-9658703-2-4; 0-8218-2695-6},
   MRCLASS = {47B80 (15A52 30E25 33D45 37K10 42C05 47B36 60F99)},
  MRNUMBER = {1677884},
MRREVIEWER = {Alexander Vladimirovich Kitaev},
}

@book {For,
    AUTHOR = {Forrester, P. J.},
     TITLE = {Log-gases and random matrices},
    SERIES = {London Mathematical Society Monographs Series},
    VOLUME = {34},
 PUBLISHER = {Princeton University Press, Princeton, NJ},
      YEAR = {2010},
     PAGES = {xiv+791},
      ISBN = {978-0-691-12829-0},
   MRCLASS = {82-02 (33C45 60B20 82B05 82B41 82B44)},
  MRNUMBER = {2641363},
MRREVIEWER = {Steven Joel Miller},
       DOI = {10.1515/9781400835416},
       URL = {https://doi.org/10.1515/9781400835416},
}

@book {HKPV,
    AUTHOR = {Hough, J. Ben and Krishnapur, M. and Peres, Y. and
              Vir\'{a}g, B.},
     TITLE = {Zeros of {G}aussian analytic functions and determinantal point processes},
    SERIES = {University Lecture Series},
    VOLUME = {51},
 PUBLISHER = {American Mathematical Society, Providence, RI},
      YEAR = {2009},
       DOI = {10.1090/ulect/051}
}

@misc{Ioo25,
      title={Partial {B}ergman kernels and determinantal point processes on {K}\"ahler manifolds}, 
      author={Louis Ioos},
      year={2025},
      eprint={2511.20539},
      note={arxiv:2511.20539},
      archivePrefix={arXiv},
      primaryClass={math.DG},
      url={https://arxiv.org/abs/2511.20539}, 
}

@incollection{Joh,
  author    = {Johansson, Kurt},
  title     = {Random matrices and determinantal processes},
  booktitle = {Mathematical Statistical Physics},
  series    = {Les Houches},
  volume    = {83},
  pages     = {1--55},
  publisher = {Elsevier},
  address   = {Amsterdam},
  year      = {2006}
}

@article {Kle,
    AUTHOR = {Klevtsov, Semyon},
     TITLE = {Random normal matrices, {B}ergman kernel and projective
              embeddings},
   JOURNAL = {J. High Energy Phys.},
  FJOURNAL = {Journal of High Energy Physics},
      YEAR = {2014},
    NUMBER = {1},
     PAGES = {133, front matter+18},
      ISSN = {1126-6708},
   MRCLASS = {60B20},
  MRNUMBER = {3599064},
       DOI = {10.1007/JHEP01(2014)133},
       URL = {https://doi.org/10.1007/JHEP01(2014)133},
}

@incollection {Kle2,
    AUTHOR = {Klevtsov, Semyon},
     TITLE = {Geometry and large {$N$} limits in {L}aughlin states},
 BOOKTITLE = {Travaux math\'{e}matiques. {V}ol. {XXIV}},
    SERIES = {Trav. Math.},
    VOLUME = {24},
     PAGES = {63--127},
 PUBLISHER = {Fac. Sci. Technol. Commun. Univ. Luxemb., Luxembourg},
      YEAR = {2016},
   MRCLASS = {81V10 (14H81 32L81)},
  MRNUMBER = {3643934},
MRREVIEWER = {Akira Asada},
}

@article {KMMW,
    AUTHOR = {Klevtsov, Semyon and Ma, Xiaonan and Marinescu, George and
              Wiegmann, Paul},
     TITLE = {Quantum {H}all effect and {Q}uillen metric},
   JOURNAL = {Comm. Math. Phys.},
  FJOURNAL = {Communications in Mathematical Physics},
    VOLUME = {349},
      YEAR = {2017},
    NUMBER = {3},
     PAGES = {819--855},
      ISSN = {0010-3616},
   MRCLASS = {81V70 (58Z05)},
  MRNUMBER = {3602817},
       DOI = {10.1007/s00220-016-2789-2},
       URL = {https://doi.org/10.1007/s00220-016-2789-2},
}

@article {Lyo,
    AUTHOR = {Lyons, Russell},
     TITLE = {Determinantal probability measures},
   JOURNAL = {Publ. Math. Inst. Hautes \'{E}tudes Sci.},
  FJOURNAL = {Publications Math\'{e}matiques. Institut de Hautes \'{E}tudes
              Scientifiques},
    NUMBER = {98},
      YEAR = {2003},
     PAGES = {167--212},
      ISSN = {0073-8301},
   MRCLASS = {60D05 (05C05 60C05 60E15)},
  MRNUMBER = {2031202},
MRREVIEWER = {Lutz Peter Klotz},
       DOI = {10.1007/s10240-003-0016-0},
       URL = {https://doi.org/10.1007/s10240-003-0016-0},
}

@book {MaMa,
    AUTHOR = {Ma, Xiaonan and Marinescu, George},
     TITLE = {Holomorphic {M}orse inequalities and {B}ergman kernels},
    SERIES = {Progress in Mathematics},
    VOLUME = {254},
 PUBLISHER = {Birkh\"{a}user Verlag, Basel},
      YEAR = {2007},
     PAGES = {xiv+422},
      ISBN = {978-3-7643-8096-0},
   MRCLASS = {32L20 (32A25 58J20 58J35 58J52 58J60)},
  MRNUMBER = {2339952},
MRREVIEWER = {David Borthwick},
}

@article {MaMa2,
    AUTHOR = {Ma, Xiaonan and Marinescu, George},
     TITLE = {Toeplitz operators on symplectic manifolds},
   JOURNAL = {J. Geom. Anal.},
  FJOURNAL = {Journal of Geometric Analysis},
    VOLUME = {18},
      YEAR = {2008},
    NUMBER = {2},
     PAGES = {565--611},
      ISSN = {1050-6926,1559-002X},
   MRCLASS = {53D50 (32A25 47B35)},
  MRNUMBER = {2393271},
MRREVIEWER = {Miroslav\ Engli\v{s}},
       DOI = {10.1007/s12220-008-9022-2},
       URL = {https://doi.org/10.1007/s12220-008-9022-2},
}

@article {Mac2,
    AUTHOR = {Macchi, O.},
     TITLE = {The coincidence approach to stochastic point processes},
   JOURNAL = {Advances in Appl. Probability},
  FJOURNAL = {Advances in Applied Probability},
    VOLUME = {7},
      YEAR = {1975},
     PAGES = {83--122},
      ISSN = {0001-8678},
   MRCLASS = {60G55},
  MRNUMBER = {380979},
MRREVIEWER = {Jan Grandell},
       DOI = {10.2307/1425855},
       URL = {https://doi.org/10.2307/1425855},
}

@book {Meh,
    AUTHOR = {Mehta, Madan Lal},
     TITLE = {Random matrices},
    SERIES = {Pure and Applied Mathematics (Amsterdam)},
    VOLUME = {142},
   EDITION = {Third},
 PUBLISHER = {Elsevier/Academic Press, Amsterdam},
      YEAR = {2004},
     PAGES = {xviii+688},
      ISBN = {0-12-088409-7},
   MRCLASS = {82-02 (15-02 15A52 60B99 60K35 82B41)},
  MRNUMBER = {2129906},
}

@article {Sos,
    AUTHOR = {Soshnikov, A.},
     TITLE = {Determinantal random point fields},
   JOURNAL = {Uspekhi Mat. Nauk},
  FJOURNAL = {Uspekhi Matematicheskikh Nauk},
    VOLUME = {55},
      YEAR = {2000},
    NUMBER = {5(335)},
     PAGES = {107--160},
      ISSN = {0042-1316},
   MRCLASS = {60G55 (60F05 60K05)},
  MRNUMBER = {1799012},
MRREVIEWER = {Boris A. Khoruzhenko},
       DOI = {10.1070/rm2000v055n05ABEH000321},
       URL = {https://doi.org/10.1070/rm2000v055n05ABEH000321},
}

@article {ShiTak,
    AUTHOR = {Shirai, Tomoyuki and Takahashi, Yoichiro},
     TITLE = {Random point fields associated with certain {F}redholm
              determinants. {I}. {F}ermion, {P}oisson and boson point
              processes},
   JOURNAL = {J. Funct. Anal.},
  FJOURNAL = {Journal of Functional Analysis},
    VOLUME = {205},
      YEAR = {2003},
    NUMBER = {2},
     PAGES = {414--463},
      ISSN = {0022-1236,1096-0783},
   MRCLASS = {60G55 (15A52 81S25 82B05 82B44)},
  MRNUMBER = {2018415},
MRREVIEWER = {Oleksiy\ Khorunzhiy},
       DOI = {10.1016/S0022-1236(03)00171-X},
       URL = {https://doi.org/10.1016/S0022-1236(03)00171-X},
}

\end{document}